\documentclass{amsart}
\usepackage{amssymb}
\usepackage{amsfonts}

\usepackage{graphicx}
\usepackage{amsmath}
\usepackage{amsthm}
\usepackage{subfigure}

\newtheorem{theorem}{Theorem}[section]
\newtheorem{proposition}[theorem]{Proposition}
\newtheorem{lemma}[theorem]{Lemma}

\theoremstyle{definition}

\newtheorem*{xrem}{Remark}
\newtheorem*{nt}{Notation}

\newtheorem{thm}{Theorem}[section]

\newtheorem{defn}[thm]{Definition}

\newtheorem*{rem*}{Remark}

\numberwithin{equation}{section}

\DeclareMathOperator {\Symp} {Symp}

\DeclareMathOperator {\Span} {span}

\def \algrest {\left [\Symp (\mathbb R^{2n})\right ]_{N}}

\def \algrestall {\bigl [\Lambda ^2(\mathbb R^{2n})\bigr ]_N}

\def \algrestclosed {\bigl [ Z ^2(\mathbb R^{2n})\bigr ]_N}

\begin{document}






\title[Symplectic $U_7, U_8, U_9$ singularities] {Symplectic $U_7, U_8$ and $U_9$ singularities}

\author{\.{Z}aneta Tr\c{e}bska}
\address{Warsaw University of Technology\newline
Faculty of Mathematics and Information Science\newline
Koszykowa 75, 00-662 Warszawa, Poland}

\email{ztrebska@mini.pw.edu.pl}

\subjclass{Primary 53D05. Secondary 14H20, 58K50, 58A10.}

\maketitle

\begin{abstract}
We use the method
of algebraic restrictions to classify symplectic $U_7, U_8$ and $U_9$
singularities. We use discrete symplectic invariants  to distinguish symplectic singularities of the curves. We also give the geometric description of symplectic classes.
\end{abstract}

\section{Introduction}
In this paper we examine the singularities which are in the list of the simple $1$-dimensional isolated complete intersection singularities in the space of dimension greater than $2$, obtained by Giusti (\cite{G}, \cite{AVG}).
Isolated complete intersection singularities (ICIS) were intensively studied by many authors (e. g. see \cite{L}), because of their interesting geometric, topological and algebraic properties.
Here using the method of algebraic restrictions we obtain the complete symplectic classification of the singularities of type $U_7, U_8$ and $U_9$. We calculate discrete symplectic invariants for symplectic orbits of the curves and we give their geometric description. It allows us to explore the specific singular nature of these classical singularities that only appears in the presence of the symplectic structure.

\smallskip

We study the symplectic classification of singular curves under the following equivalence:

\begin{defn} \label{symplecto}
Let $N_1, N_2$ be germs of subsets of symplectic space $(\mathbb{R}^{2n}, \omega)$. $N_1, N_2$ are \emph{symplectically equivalent} if there exists a symplectomorphism-germ $\Phi:(\mathbb{R}^{2n}, \omega) \rightarrow(\mathbb{R}^{2n}, \omega)$ such that $\Phi(N_1)=N_2$.
\end{defn}

We recall that  $\omega$ is a symplectic form if $\omega$ is a
smooth nondegenerate closed 2-form, and $\Phi:\mathbb{R}^{2n}
\rightarrow\mathbb{R}^{2n}$ is a symplectomorphism if $\Phi$ is
diffeomorphism and $\Phi ^* \omega=\omega$.

\smallskip

Symplectic classification of curves was initiated  
by V. I. Arnold. In \cite{Ar1} and \cite{Ar2} the author studied singular curves in symplectic and contact spaces and introduced the local symplectic and contact algebra. He 
discovered new symplectic invariants of singular curves. He proved that the $A_{2k}$
singularity of a planar curve (the orbit with respect to standard
$\mathcal A$-equivalence of parameterized curves) split into
exactly $2k+1$ symplectic singularities (orbits with respect to
symplectic equivalence of parameterized curves). He distinguished different symplectic singularities by different orders of tangency of the parameterized curve to the \emph{nearest} smooth Lagrangian submanifold. Arnold posed a problem of expressing these invariants in terms of the local algebra's interaction with the symplectic structure and he proposed calling this interaction the local symplectic algebra.

\smallskip

In \cite{IJ1} G. Ishikawa and S. Janeczko classified symplectic singularities of curves in the $2$-dimensional symplectic space. All simple curves in this classification are quasi-homogeneous.

\smallskip

We recall that a subset $N$ of $\mathbb R^m$ is \emph{quasi-homogeneous} if there exist a coordinate system $(x_1,\cdots,x_m)$ on $\mathbb R^m$ and
positive numbers $w_1,\cdots,w_m$ (called weights) such that for any
point $(y_1,\cdots,y_m)\in \mathbb R^m$ and any $t>0$
if $(y_1,\cdots,y_m)$ belongs to $N$ then the point
$(t^{w_1}y_1,\cdots,t^{w_m}y_m)$ belongs to $N$.

\smallskip

The generalization of results in \cite{IJ1} to volume-preserving classification of singular varieties and maps  in arbitrary dimensions was obtained in \cite{DR}. A symplectic form on a $2$-dimensional manifold is a special case of a volume form on a smooth manifold.

The stably simple symplectic singularities of parameterized curves (in the $\mathbb C$-analytic category) were studied  by P. A. Kolgushkin in \cite{K}. 


In \cite{Zh} was developed the local contact algebra. The main results were based on the notion of the algebraic restriction of a contact structure to a subset $N$ of a contact manifold.

In \cite{DJZ2} new symplectic invariants of singular quasi-homogeneous  subsets of a symplectic space were explained by the algebraic restrictions of the symplectic form to these subsets.

\smallskip

The algebraic restriction is an equivalence class of the following relation on the space of differential $k$-forms:

Differential $k$-forms $\omega_1$ and $\omega_2$ have the same
\emph{algebraic restriction} to a subset $N$ if
$\omega_1-\omega_2=\alpha+d\beta$, where $\alpha$ is a $k$-form
vanishing on $N$ and $\beta$ is a $(k-1)$-form vanishing on $N$.

\smallskip

In \cite{DJZ2} the generalization of Darboux-Givental theorem (\cite{ArGi})
to germs of arbitrary subsets of the symplectic space was obtained. This result reduces
the problem of symplectic classification of germs of quasi-homo\-ge\-neous subsets to
the problem of classification of algebraic restrictions of symplectic
forms to these subsets. For non-quasi-homogeneous subsets there is one more cohomological invariant apart from the algebraic restriction (\cite{DJZ2}, \cite{DJZ1}). The dimension of the space of algebraic restrictions of closed $2$-forms to a $1$-dimensional quasi-homogeneous isolated complete
intersection singularity $C$ is equal to the multiplicity of $C$ (\cite{DJZ2}). In \cite{D1} it was proved that the space of algebraic restrictions of closed $2$-forms to a $1$-dimensional (singular) analytic variety is finite-dimensional.
In \cite{DJZ2} the method of algebraic restrictions was applied to various classification problems in a symplectic space. In particular the complete symplectic classification of classical $A$-$D$-$E$ singularities of planar curves and $S_5$ singularity was obtained. Most of different symplectic singularity classes were distinguished by new discrete symplectic invariants: the index of isotropy and the symplectic multiplicity.

\smallskip

In \cite{DT1} following ideas from \cite{Ar1} and \cite{D1}  new
discrete symplectic invariants - the Lagrangian tangency orders
were introduced and used to distinguish symplectic
singularities of simple planar curves of type $A$-$D$-$E$, symplectic
 $T_7$ and $T_8$ singularities.

 The complete symplectic classification of
the isolated complete intersection singularities $S_{\mu}$
for $\mu>5$ and $W_8$, $W_9$ singularities were obtained in \cite{DT2} and \cite{Tr} respectively.

 \smallskip

 The method of algebraic restrictions was successfully used by W. Domitrz in \cite{D2} to classify the 0-dimensional ICIS 
  (multiple points) in a symplectic space.

 \smallskip

  In this paper we obtain the detailed symplectic classification of the  $U_7, U_8$ and the $U_9$ singularities.
The paper is organized as follows.  In Section \ref{discrete} we recall  discrete symplectic invariants (the symplectic multiplicity, the index of isotropy and the Lagrangian tangency orders).
Symplectic classification of the $U_7, U_8$ and the $U_9$ singularity is presented in Sections \ref{sec-u7}, \ref{sec-u8} and \ref{sec-u9} respectively. The symplectic sub-orbits of this singularities are listed in Theorems \ref{u7-main}, \ref{u8-main} and \ref{u9-main}. Discrete symplectic invariants for the symplectic classes are calculated in Theorems \ref{lagr-u7}, \ref{lagr-u8} and \ref{lagr-u9}. The geometric descriptions of the symplectic orbits are presented in Theorems \ref{geom-cond-u8}, \ref{geom-cond-u8} and \ref{geom-cond-u9}.
In Section \ref{proofs} we recall the method of algebraic restrictions and use it to classify  symplectic singularities.

\section{Discrete symplectic invariants}\label{discrete}

We can use discrete symplectic invariants to characterize
symplectic singularity classes.

 The first invariant is a symplectic
multiplicity (\cite{DJZ2}) introduced  in \cite{IJ1} as a
symplectic defect of a curve.

\medskip

Let $N$ be a germ of a subvariety of $(\mathbb R^{2n},\omega)$.

\begin{defn}
\label{def-mu}
 The \emph{symplectic multiplicity}, $\mu^{sym}(N)$ of  $N$ is the codimension of
 the symplectic orbit of $N$ in the orbit of $N$ with respect to the action of the group of local diffeomorphisms.
\end{defn}

The second invariant is the index of isotropy \cite{DJZ2}.

\begin{defn}
The \emph{index of isotropy}, $ind(N)$ of $N$ is the maximal
order of vanishing of the $2$-forms $\omega \vert _{TM}$ over all
smooth submanifolds $M$ containing $N$.
\end{defn}

This invariant has geometrical interpretation. An equivalent definition is as follows: the index of isotropy of $N$ is the maximal order of tangency between non-singular submanifolds containing $N$ and non-singular isotropic submanifolds of the same dimension.  The index of isotropy is equal to $0$ if $N$ is not contained in any non-singular submanifold which is tangent to some isotropic submanifold of the same dimension. If $N$ is contained in a non-singular Lagrangian submanifold then the index of isotropy is $\infty $.


The symplectic multiplicity and the index of isotropy can be described in terms of algebraic restrictions (Propositions \ref{sm} and \ref{ii} in Section \ref{proofs}).

\medskip

There is one more discrete symplectic invariant, introduced in \cite{D1} (following ideas from \cite{Ar2}) which is defined specifically for a parameterized curve. This is the maximal
tangency order of a curve $f:\mathbb R\rightarrow M$ to a smooth Lagrangian submanifold. If $H_1=...=H_n=0$ define a smooth submanifold $L$ in the symplectic space then the tangency order of
a curve $f:\mathbb R\rightarrow M$ to $L$ is the minimum of the orders of vanishing at $0$ of functions $H_1\circ f,\cdots, H_n\circ f$. We denote the tangency order of $f$ to $L$ by $t(f,L)$.

\begin{defn}
The \emph{Lagrangian tangency order} $Lt(f)$ \emph{ of a curve} $f$ is the
maximum of $t(f,L)$ over all smooth Lagrangian submanifolds $L$ of
the symplectic space.
\end{defn}

The Lagrangian tangency order of the quasi-homogeneous curve in a symplectic space can also be  expressed in terms of algebraic restrictions  (Proposition \ref{lto} in Section \ref{proofs}).

\medskip

In \cite{DT1}  the above invariant was generalized for germs of
curves and multi-germs of curves which may be parameterized
analytically since the Lagrangian tangency order is the same for every
'good' analytic parametrization of a curve.

\medskip

Consider a multi-germ $(f_i)_{i\in\{1,\cdots,r\}}$ of analytically
parameterized curves $f_i$.  We have $r$-tuples $(t(f_1,L), \cdots,
t(f_r,L))$ for any smooth submanifold $L$ in the
symplectic space.

\begin{defn}
For any $I\subseteq \{1,\cdots, r\}$ we define \emph{the tangency order of the multi-germ } $(f_i)_{i\in I}$ to $L$:
$$t[(f_i)_{i\in\ I},L]=\min_{i\in\ I} t(f_i,L).$$
\end{defn}

\begin{defn}
The \emph{Lagrangian tangency order} $Lt((f_i)_{i\in\ I})$ \emph{of a multi-germ } $(f_i)_{i\in I}$ is the maximum of $t[(f_i)_{i\in\ I},L]$ over all smooth Lagrangian submanifolds $L$ of the symplectic space.
\end{defn}



\section{Symplectic $U_7$-singularities}\label{sec-u7}

 Denote by $(U_7)$  the class of varieties in  a fixed symplectic space $(\mathbb R^{2n}, \omega )$ which are diffeomorphic to
 \begin{equation}
\label{defu7} U_7=\{x\in \mathbb R ^{2n\geq
4}\,:x_1^2+x_2 x_3=x_1 x_2+x_3^3=x_{\geq 4}=0\}.\end{equation}

\noindent This is the simple  $1$-dimensional isolated complete intersection singularity $U_7$ (\cite{G}, \cite{AVG}).
Here $N$ is quasi-homogeneous with weights $w(x_1)=4,\; w(x_2)=5$,\; $w(x_3)=3$.

\medskip

We used the method of algebraic restrictions to obtain the complete classification of symplectic singularities of $(U_7)$ presented in the following theorem.

\begin{theorem}\label{u7-main}
Any submanifold of the symplectic space $(\mathbb R^{2n},\sum_{i=1}^n dp_i \wedge dq_i)$ where $n\geq3$ (respectively $n=2$) which is diffeomorphic to $U_7$ is symplectically equivalent to one and only one of the normal forms $U_7^i, i = 0,1,\cdots ,7$  (respectively $i=0,1,2)$
  listed below. The parameters $c, c_1, c_2$ of the normal forms are moduli:

\smallskip
\begin{small}
\noindent $U_7^0$: \ $p_1^2 + p_2q_1 = 0, \ \ p_1p_2 +q_1^3= 0, \
\ q_2 = c_1q_1 + c_2p_1, \ \ p_{\ge 3} = q_{\ge 3} = 0$;

\smallskip

\noindent $U_7^1$: \ $p_2^2 \pm p_1q_1 = 0, \ \ p_1p_2 \pm q_1^3= 0, \
\ q_2 = c_1p_1 + \frac{c_2}{2}q_1^2, \ \ p_{\ge 3} = q_{\ge 3} = 0$;

\smallskip

\noindent $U_7^2$: \ $p_1^2 + q_1q_2 = 0, \ \ p_1q_1+q_2^3 = 0, \
\ p_2 = c_1p_1q_2+\frac{c_2}{2}p_1^2, \ \ p_{\ge 3} = q_{\ge 3} = 0$;

\smallskip

\noindent $U_7^3$: \ $p_1^2 + p_2p_3 = 0,  \ p_1p_2 +p_3^3 = 0, \ q_1 = cp_1p_3, q_2=0, q_3=\pm p_1p_3,  \ p_{\geq 4} = q_{\geq 4} = 0$;

\smallskip

\noindent $U_7^4$: \ $p_1^2 + p_2p_3 = 0,  \ p_1p_2 +p_3^3 = 0, \ q_1 = \frac{c}{3}p_3^3, q_2=0, q_3=-\frac{1}{2}p_1^2,  \ p_{\geq 4} = q_{\geq 4} = 0$;

\smallskip

\noindent $U_7^5$: \ $p_1^2 + p_2p_3 = 0,  \ p_1p_2 +p_3^3 = 0, \ q_1 = -\frac{c}{2}p_1p_3^2, q_2=0, q_3=-p_1p_3^2,  \ p_{\geq 4} = q_{\geq 4} = 0$;

\smallskip

\noindent $U_7^6$: \ $p_1^2 + p_2p_3 = 0,  \ p_1p_2 +p_3^3 = 0, \ q_1 =0, q_2=0, q_3= \mp\frac{1}{2}p_1^2p_3,   \ p_{\ge 4} = q_{\ge 4} = 0;$

\smallskip

\noindent $U_7^7$: \ $p_1^2 + p_2p_3 = 0,  \ p_1p_2 +p_3^3 = 0, \
q_{\ge 1} = p_{\ge 4} = 0.$
\end{small}
\end{theorem}


\subsection{Distinguishing symplectic classes of $U_7$ by the Lagrangian tangency orders}
\label{u7-lagr}

A curve $N\in (U_7)$ may be described as a union of two parametrical branches $B_1$ and $B_2$. The branch $B_1$ is smooth so it is contained in some Lagrangian submanifold and thus $Lt(B_1)=\infty$. The branch $B_2$ is singular. The parametrizations of branches are given in Table \ref{tabu7-lagr}. To characterize the symplectic classes  we use the following invariants:
\begin{itemize}
  \item $Lt=Lt(B_1,B_2)=\max\limits _\mathcal{L} (\min \{t(B_1,\mathcal{L}),t(B_2,\mathcal{L})\}),$

  \item $L_{2}=Lt(B_2)=\max\limits _\mathcal{L}\, t(B_2,\mathcal{L}),$
\end{itemize}
Here $L$ is a smooth Lagrangian submanifold of the symplectic space.
We  also compute $ind$ (the index of isotropy of $N$) and $ind_2$ (the index of isotropy of the singular component).
\begin{theorem}
\label{lagr-u7} Any stratified submanifold $N\in (U_7)$ of a symplectic space $(\mathbb R^{2n}, \omega_0)$ with the canonical coordinates $(p_1, q_1, \cdots, p_n, q_n)$ is symplectically equivalent to one and only one of the curves presented in the second column of Table \ref{tabu7-lagr}. 
The  indices of isotropy and the Lagrangian tangency orders of the curve $N$   are presented  in the third, fourth, fifth and sixth column of Table \ref{tabu7-lagr}.
\end{theorem}
\renewcommand*{\arraystretch}{1.3}
\begin{center}
\begin{table}[h]
    \begin{small}
    \noindent
    \begin{tabular}{|p{1.2cm}|p{5.8cm} r|c|c|c|c|}
                     \hline
    class &  parametrization of branches of $N$ & & $ind$ & $ind_2$ & $Lt$   &  $L_{2}$  \\ \hline

 $(U_7)^0$ & $B_1:(0,0,t,0,0,\cdots )$ & if $c_1\ne 0$ & $0$ & $0$ & $3$ & $4$  \\ \cline{3-7}
     $2n\ge 4$ & $B_2:(t^4,-t^3,t^5,-c_1t^3-c_2t^4,0,\cdots )$ & if $c_1=0$ & $0$ & $0$ & $4$ & $4$  \\ \hline

     $(U_7)^1$ & $B_1:(t,0,0,c_1t,,0,\cdots )$ & & $0$ & $0$ & $3$ & $5$  \\
     $2n\ge 4$ & $B_2:(t^5,\mp t^3,t^4,c_1t^5+\frac{c_2}{2}t^6,0,\cdots )$ & & &  &  &   \\ \hline

   $(U_7)^2$  & $B_1:(0,t,0,0,0,\cdots )$ & & $0$ & $0$ & $4$ & $5$    \\
    $2n\ge 4$  & $B_2:(t^4,t^5,-c_1t^7+\frac{c_2}{2}t^8,-t^3,0,\cdots )$ & & &  &  &   \\ \hline

    $(U_7)^3$ & $B_1:(0,0,t,0,0,0,0,\cdots )$ & & $1$ & $1$ & $7$  & $7$ \\
    $2n\ge 6$ & $B_2:(t^4,-ct^7,t^5,0,-t^3,\pm t^7,0,\cdots )$ &  & &   &   &  \\ \hline

    $(U_7)^4$ & $B_1:(0,0,t,0,0,0,0,\cdots )$ & & $1$ & $1$ & $8$  & $8$ \\
    $2n\ge 6$ & $B_2:(t^4,-\frac{c}{3}t^9,t^5,0,-t^3,-\frac{1}{2} t^{8},0,\cdots )$ &  & &   &   &  \\ \hline

    $(U_7)^5$ & $B_1:(0,0,t,0,0,0,0,\cdots )$ & & $2$ & $\infty$ & $10$  & $\infty$ \\
    $2n\ge 6$ & $B_2:(t^4,-\frac{c}{2}t^{10},t^5,0,-t^3,-t^{10},0,\cdots )$ & &  &   &   &  \\ \hline

     $(U_7)^6$ & $B_1:(0,0,t,0,0,0,0,\cdots )$ & & $2$ & $\infty$ & $11$  & $\infty$ \\
    $2n\ge 6$ & $B_2:(t^4,0,t^5,0,-t^3,\pm \frac{1}{2}t^{11},0,\cdots )$ & &  &   &   &  \\ \hline

    $(U_7)^7$ & $B_1:(0,0,t,0,0,0,0,\cdots )$ & & $\infty$ & $\infty$ & $\infty$  & $\infty$ \\
    $2n\ge 6$ & $B_2:(t^4,0,t^5,0,-t^3,0,0,\cdots )$ &  & &   &   &  \\ \hline
\end{tabular}

\medskip

\caption{\small The symplectic invariants for symplectic classes of $U_7$ singularity.}\label{tabu7-lagr}

\end{small}
\end{table}
\end{center}



\begin{xrem}
The comparison of invariants presented in Table \ref{tabu7-lagr} shows  that the Lagrangian tangency orders distinguish more symplectic classes than the respective indices of isotropy.

 The most of invariants can be calculated by knowing algebraic restrictions for the symplectic classes. We use Proposition \ref{ii} to calculate the indices of isotropy. $L_{2}$ is calculated by using Proposition \ref{lto} for the singular branch. $Lt$ is computed by  applying directly the definition of the Lagrangian tangency order and finding a Lagrangian submanifold the nearest to the curve $N$.
\end{xrem}

\smallskip

\subsection{Identifying the classes $(U_7)^i$ by geometric conditions}
\label{u7-geom_cond} $ $

\noindent We can characterize the symplectic classes $(U_7)^i$ by geometric conditions independent of any local coordinate system.

Let $N\in (U_7)$. Denote by $W$ the tangent space at $0$ to some 
non-singular \linebreak $3$-manifold containing $N$. We can define the following subspaces of this space:

  $\ell_1$ -- the tangent line at $0$ to the nonsingular branch $B_1$,

  $\ell_2$ -- the tangent line at $0$ to the singular branch $B_2$

  $V$ -- the $2$-space tangent at $0$ to the singular branch $B_2$.

  For $N=U_7=$(\ref{defu7}) it is easy to calculate that
$W\!=\!\Span (\partial /\partial x_1,\partial /\partial x_2,\partial /\partial x_3)$,  \newline and
$\ell\!_1=\!\Span (\partial /\partial x_2), \ \ell\!_2=\!\Span (\partial /\partial x_3), \
V\!=\!\Span (\partial /\partial x_1,\partial /\partial x_3)$.

\smallskip

The classes $(U_7)^i$ satisfy special conditions in terms of the restriction $\omega\vert_ W $, where $\omega$ is the symplectic form.

\begin{theorem}
\label{geom-cond-u7} Any stratified submanifold $N\in (U_7)$ of a
symplectic space $(\mathbb R^{2n}, \omega)$ belongs to the class
$(U_7)^i$ if and only if the couple $(N, \omega)$ satisfies the corresponding
conditions in the last column of Table \ref{tabu7-geom}.
\end{theorem}

\begin{center}
\begin{table}[h]
    \noindent
    \begin{tabular}{|p{1.1cm}|p{5.6cm}|p{4.7cm}|}
              \hline
    class &  normal form & geometric conditions  \\ \hline

  $(U_7)^0$   &  $[U_7]^0_0: [\theta _1 + c_1\theta _2 + c_2\theta _3]_{U_7},\; c_1\ne 0$
    &  $\omega|_V \neq 0$ and $\omega|_{\ell_1+\ell_2} \neq 0$ \\ \cline{2-3}
    &  $[U_7]^0_1: [\theta _1 + c_2\theta _3]_{U_7}$
    &  $\omega|_V \neq 0$ and $\omega|_{\ell_1+\ell_2} = 0$ \\
     \hline
  $(U_7)^1$ & $[U_7]^1: [\pm\theta _2 + c_1\theta _3 + c_2\theta _4]_{U_7}$ & $\omega|_V=0$ but \;$\ker\omega \ne \ell_2$    \\ \hline

 $(U_7)^{2}$ & $[U_7]^2: [\theta _3 + c_1\theta_4+c_2\theta_5]_{U_7}$ & $\omega|_V=0$ and $\ker\omega = \ell_2$ \\ \hline

 & & $\omega\vert_ W = 0$\\ \hline
$(U_7)^3$  & $[U_7]^3: [\pm\theta _4 +c\theta_5]_{U_7}$  &  $Lt=L_2=7$  \\ \hline
$(U_7)^4$ & $[U_7]^4: [\theta _5 + c\theta _6]_{U_7}$   &  $Lt=L_2=8$   \\ \hline
$(U_7)^5$ & $[U_7]^5: [\theta _6 + c\theta _7]_{U_7}$  &  $Lt=10,\;L_2=\infty$   \\ \hline
$(U_7)^6$ & $[U_7]^6: [\pm\theta _7]_{U_7}$  &  $Lt=11,\; L_2=\infty$   \\ \hline
$(U_7)^7$ & $[U_7]^7: [0]_{U_7}$ &  $N$ is contained in a smooth \newline Lagrangian submanifold \\ \hline
    \end{tabular}

\medskip

\caption{\small Geometric interpretation of singularity classes of $U_7$. ($W$ is the tangent space to a non-singular 3-dimensional manifold in $(\mathbb R^{2n\geq4}, \omega)$ containing $N\in(U_7)$. The forms $\theta_1,\ldots,\theta_7$ are described in  Theorem ~\ref{u7-baza} on the page \pageref{u7-baza}.)}\label{tabu7-geom}
\end{table}
\end{center}

\begin{proof}[Sketch of the proof of Theorem \ref{geom-cond-u7}]
 We have to show that the conditions in the row of $(U_7)^i$ are satisfied for any $N\in (U_7)^i$.
 Each of the conditions in the last column of Table \ref{tabu7-geom}
is invariant with respect to the action of the group of diffeomorphisms in the space of pairs $(N,\omega)$.
 Because each of these conditions depends only on the algebraic restriction $[\omega ]_N$ we can take the simplest $2$-forms $\omega ^i$ representing
the normal forms $[U_7]^i$ for algebraic restrictions and we can check that the pair $(U_7,\omega=\omega ^i)$ satisfies the condition
in the last column of Table \ref{tabu7-geom}. By simple calculation and observation of the Lagrangian tangency orders we obtain that the conditions corresponding to the classes $(U_7)^i$ are satisfied.
\end{proof}


\section{Symplectic $U_8$-singularities}\label{sec-u8}

 Denote by $(U_8)$  the class of varieties in  a fixed symplectic space $(\mathbb R^{2n}, \omega )$ which are diffeomorphic to

 \begin{equation}
\label{defu8} U_8=\{x\in \mathbb R ^{2n\geq
4}\,:x_1^2+x_2 x_3=x_1 x_2+x_1x_3^2=x_{\geq 4}=0\}.\end{equation}

\noindent This is the simple  $1$-dimensional isolated complete intersection singularity $U_8$ (\cite{G}, \cite{AVG}).
Here $N$ is quasi-homogeneous with weights $w(x_1)=3,\; w(x_2)=4$,\; $w(x_3)=2$.

\medskip

We used the method of algebraic restrictions to obtain the complete classification of symplectic singularities of $(U_8)$ presented in the following theorem.

\begin{theorem}\label{u8-main}
Any submanifold of the symplectic space $(\mathbb R^{2n},\sum_{i=1}^n dp_i \wedge dq_i)$ where $n\geq3$ (respectively $n=2$) which is diffeomorphic to $U_8$ is symplectically equivalent to one and only one of the normal forms $U_8^i, i = 0,1,\cdots ,8$ 
  listed below. The~parameters $c, c_1, c_2$ of the normal forms are moduli:

\medskip
\smallskip
\begin{small}

\noindent $U_8^0$: \ $p_1^2 + p_2q_1 = 0, \ \ p_1p_2 +p_1q_1^3= 0, \
\ q_2 = c_1q_1 - c_2p_1, \ \ p_{\ge 3} = q_{\ge 3} = 0$;

\smallskip

\noindent $U_8^1$: \ $p_1^2 \pm p_2q_2 = 0, \ \ p_1p_2 + p_1q_2^2= 0, \
\ q_1 = c_1p_2 + \frac{c_2}{2}q_2^2, \ \ p_{\ge 3} = q_{\ge 3} = 0$;

\smallskip

\noindent $U_8^2$: \ $p_1^2 + q_1q_2 = 0, \ \ p_1q_1+p_1q_2^2 = 0, \
\ p_2 = c_1p_1q_2+\frac{c_2}{2}p_1^2, \ \ p_{\ge 3} = q_{\ge 3} = 0$;

\smallskip

\noindent $U_8$$^{3,0}_5$: \ $p_1^2 + q_1q_2 = 0, \ \ p_1q_1+p_1q_2^2 = 0, \
\ p_2 = -\frac{1}{3}p_1q_2+  \frac{c_1}{2}p_1^2+c_2p_1q_2^2, \ \ p_{\ge 3} = q_{\ge 3} = 0$;

\smallskip

\noindent $U_8$$^{3,0}_{\infty}$: \ $p_1^2 + q_1q_2 = 0, \ \ p_1q_1+p_1q_2^2 = 0, \
\ p_2 = 2p_1q_2+  \frac{c_1}{2}p_1^2+\frac{c_2}{2}p_1^2q_2, \ \ p_{\ge 3} = q_{\ge 3} = 0$;

\smallskip

\noindent $U_8$$^{3,1}$: \ $p_1^2 + p_2p_3 = 0,  \ p_1p_2 +p_1p_3^2 = 0, \ q_1 =  q_2=0, q_3=-p_1p_3-\frac{c}{2}p_1^2,  \ p_{>3} = q_{>3} = 0$;

\smallskip

\noindent $U_8^4$: \ $p_1^2 + p_2p_3 = 0,  \ p_1p_2 +p_1p_3^2 = 0, \ q_1 =  q_2=0, q_3=\mp \frac{1}{2}p_1^2-cp_1p_3^2,  \ p_{>3} = q_{>3} = 0$;

\smallskip

\noindent $U_8^5$: \ $p_1^2 + p_2p_3 = 0,  \ p_1p_2 +p_1p_3^2 = 0, \ q_1 = q_2=0, q_3=-p_1p_3^2-\frac{c}{2}p_1^2p_3,  \ p_{>3} = q_{>3} = 0$;

\smallskip

\noindent $U_8^6$: \ $p_1^2 + p_2p_3 = 0,  \ p_1p_2 +p_1p_3^2 = 0, \ q_1 = q_2=0, q_3=\mp \frac{1}{2}p_1^2p_3+cp_1p_3^3,  \ p_{>3} = q_{>3} = 0$;

\smallskip

\noindent $U_8^7$: \ $p_1^2 + p_2p_3 = 0,  \ p_1p_2 +p_1p_3^2 = 0, \  q_1 = q_2=0, q_3 =-p_1p_3^3, \ p_{>3} = q_{>3} = 0$;,

\smallskip

\noindent $U_8^{8}$: \ $p_1^2 + p_2p_3 = 0,  \ p_1p_2 +p_1p_3^2 = 0, \
q_{\ge 1} = p_{\ge 4} = 0.$
\end{small}
\end{theorem}
\medskip


\subsection{Distinguishing symplectic classes of $U_8$ by the Lagrangian tangency orders}
\label{u8-lagr}

A curve $N\in (U_8)$ may be described as a union of three parametrical branches $B_1, B_2$ and $B_3$. Branches $B_1, B_2$ are smooth and their union is an invariant component diffeomorphic to $A_1$ singularity and the branch $B_3$ is diffeomorphic to $A_2$ singularity. Their parametrizations are given in Table \ref{tabu8-lagr}. To characterize the symplectic classes  we use the following invariants:
\begin{itemize}
  \item $Lt=Lt(B_1,B_2,B_3)=\max\limits _\mathcal{L} (\min \{t(B_1,\mathcal{L}),t(B_2,\mathcal{L}),t(B_3,\mathcal{L})\}),$
  \item $L_{1,2}=Lt(B_1,B_2)=\max\limits _\mathcal{L} (\min \{t(B_1,\mathcal{L}),t(B_2,\mathcal{L})\}),$
  \item $L_{3}=Lt(B_3)=\max\limits _\mathcal{L}\, t(B_3,\mathcal{L}),$
\end{itemize}
Here $L$ is a smooth Lagrangian submanifold of the symplectic space.
\begin{theorem}
\label{lagr-u8} Any stratified submanifold $N\in (U_8)$ of a symplectic space $(\mathbb R^{2n}, \omega_0)$ with the canonical coordinates $(p_1, q_1, \cdots, p_n, q_n)$ is symplectically equivalent to one and only one of the curves presented in the second column of Table \ref{tabu8-lagr}. 
The  index of isotropy of the curve $N$ and the Lagrangian tangency orders  are presented  in the third and fourth, fifth and sixth column of Table \ref{tabu8-lagr}.
\end{theorem}
\setlength{\tabcolsep}{1.2mm}
\renewcommand*{\arraystretch}{1.3}
\begin{center}
\begin{table}[h]

    \begin{small}
    \noindent
    \begin{tabular}{|p{1.0cm}|p{6.8cm} r|c|c|c|c|}
                      \hline
    class &  parametrization of branches of $N$  &   & $ind$ & $Lt$   &  $L_{1,2}$ & $L_3$ \\ \hline

 $(U_8)^0$ & $B_1:(0,0,t,0,0,\cdots ),\; B_2:(0,t,0,c_1t,0,\cdots)$ &  $c_1\ne 0$ & $0$ & $1$ & $1$ & $3$ \\ \cline{3-7}
     $2n\ge 4$ & $B_3:(t^3,t^2,-t^4,c_1t^2-c_2t^3,0,\cdots )$ &  $c_1=0$ & $0$ & $3$ & $\infty$ & $3$ \\ \hline

$(U_8)^1$ & $B_1\!:(0,c_1t,t,0,0,\cdots),\; B_2\!:(0,\frac{c_2}{2}t^2,0,\pm t,0,\cdots)$ &  $c_2\!\ne 2c_1$ & $0$ & $1$ & $1$ & $5$ \\ \cline{3-7}
     $2n\ge 4$ & $B_3:(t^3,(\frac{c_2}{2}-c_1)t^4,-t^4,\pm t^2,0,\cdots )$ &  $c_2\!=2c_1$ & $0$ & $1$ & $1$ & $\infty$ \\ \hline

$(U_8)^2$ & $B_1:(0,t,0,0,0,\cdots ),\; B_2:(0,0,0,t,0,\cdots)$ &  $c_1\ne 2,$ & $0$ & $3$ & $\infty$ & $5$ \\
     $2n\ge 4$ & $B_3:(t^3,-t^4,c_1t^5+\frac{c_2}{2}t^6, t^2,0,\cdots )$ &  $c_1\!\ne\!-\frac{1}{3}$ & & &  &  \\ \hline

$(U_8)^{3,0}_5$ & $B_1:(0,t,0,0,0,\cdots ),\; B_2:(0,0,0,t,0,\cdots)$ &   & $0$ & $3$ & $\infty$ & $5$ \\
     $2n\ge 4$ & $B_3:(t^3,-t^4,-\frac{1}{3}t^5+\frac{c_1}{2}t^6+c_2t^7, t^2,0,\cdots )$ &   & & &  &  \\ \hline

$(U_8)^{3,0}_{\infty}$ & $B_1:(0,t,0,0,0,\cdots ),\; B_2:(0,0,0,t,0,\cdots)$ &   & $0$ & $3$ & $\infty$ & $\infty$ \\
     $2n\ge 4$ & $B_3:(t^3,-t^4,2t^5+\frac{c_1}{2}t^6+\frac{c_2}{2}t^8, t^2,0,\cdots )$ &   & & &  &  \\ \hline

$(U_8)^{3,1}$ & $B_1:(0,0,t,0,0,0,\cdots ),\; B_2:(0,0,0,0,t,0,\cdots)$ &   & $1$ & $5$ & $\infty$ & $5$ \\
     $2n\ge 6$ & $B_3:(t^3,0,-t^4,0, t^2,-t^5-\frac{c}{2}t^6,0,\cdots )$ &   & & &  &  \\ \hline

$(U_8)^4$ & $B_1:(0,0,t,0,0,0,\cdots ),\; B_2:(0,0,0,0,t,0,\cdots)$ &   & $1$ & $6$ & $\infty$ & $\infty$ \\
     $2n\ge 6$ & $B_3:(t^3,\pm t^5,-t^4,0, t^2,-ct^7,0,\cdots )$ &   & & &  &  \\ \hline

$(U_8)^5$ & $B_1:(0,0,t,0,0,0,\cdots ),\; B_2:(0,0,0,0,t,0,\cdots)$ &   & $2$ & $7$ & $\infty$ & $\infty$ \\
     $2n\ge 6$ & $B_3:(t^3,0,-t^4,0, t^2,-t^7-\frac{c}{2}t^8,0,\cdots )$ &   & & &  &  \\ \hline

$(U_8)^6$ & $B_1:(0,0,t,0,0,0,\cdots ),\; B_2:(0,0,0,0,t,0,\cdots)$ &   & $2$ & $8$ & $\infty$ & $\infty$ \\
     $2n\ge 6$ & $B_3:(t^3,0,-t^4,0, t^2,\mp \frac{1}{2}t^8+ct^9,0,\cdots )$ &   & & &  &  \\ \hline

$(U_8)^7$ & $B_1:(0,0,t,0,0,0,\cdots ),\; B_2:(0,0,0,0,t,0,\cdots)$ &   & $3$ & $9$ & $\infty$ & $\infty$ \\
     $2n\ge 6$ & $B_3:(t^3,0,-t^4,0, t^2,-t^9,0,\cdots )$ &   & & &  &  \\ \hline

$(U_8)^{8}$ & $B_1:(0,0,t,0,0,0,\cdots ),\; B_2:(0,0,0,0,t,0,\cdots)$ &   & $\infty$ & $\infty$ & $\infty$ & $\infty$ \\
     $2n\ge 6$ & $B_3:(t^3,0,-t^4,0, t^2,0,0,\cdots )$ &   & & &  &  \\ \hline

\end{tabular}

\smallskip

\caption{\small The symplectic invariants for symplectic classes of $U_8$ singularity.}\label{tabu8-lagr}

\end{small}
\end{table}
\end{center}
\begin{xrem}
The comparison of invariants presented in Table \ref{tabu8-lagr} shows  that the Lagrangian tangency order distinguishes more symplectic classes than the index of isotropy. Symplectic classes $(U_8)^{2}$ and $(U_8)^{3,0}_5$ can  be distinguished by the symplectic multiplicity.

 The invariants can be calculated by knowing algebraic restrictions for the symplectic classes. We use Proposition \ref{ii} to calculate the index of isotropy. The invariants $L_{1,2}$ and $L_{3}$ we can calculate knowing the respective Lagrangian tangency orders for $A_1$ and $A_2$ singularities. $Lt$ is computed by  applying directly the definition of the Lagrangian tangency order and finding a Lagrangian submanifold the nearest to the curve $N$.
 \end{xrem}


\subsection{Geometric conditions for the classes $(U_8)^i$}
\label{u8-geom_cond} $ $

\noindent We can characterize the symplectic classes $(U_8)^i$ by geometric conditions independent of any local coordinate system.

Let $N\in (U_8)$. Denote by $W$ the tangent space at $0$ to some 
non-singular \linebreak $3$-manifold containing $N$. We can define the following subspaces of this space:

  $\ell_1$ -- the tangent line at $0$ to the nonsingular branch $B_1$,

  $\ell_2$ -- the tangent line at $0$ to the  nonsingular branch $B_2$ (this line is also tangent at $0$ to the singular branch $B_3$),

  $V$ -- the $2$-space tangent at $0$ to the singular branch $B_3$.

  For $N=U_8=$(\ref{defu8}) it is easy to calculate that
$W\!=\!\Span (\partial /\partial x_1,\partial /\partial x_2,\partial /\partial x_3)$,  \newline and
$\ell\!_1=\!\Span (\partial /\partial x_2), \ \ell\!_2=\!\Span (\partial /\partial x_3), \
V\!=\!\Span (\partial /\partial x_1,\partial /\partial x_3)$.

\smallskip

The classes $(U_8)^i$ satisfy special conditions in terms of the restriction $\omega\vert_ W $, where $\omega$ is the symplectic form.

\begin{theorem}
\label{geom-cond-u8}If a stratified submanifold $N\in (U_8)$ of a
symplectic space $(\mathbb R^{2n}, \omega)$ belongs to the class
$(U_8)^i$ then the couple $(N, \omega)$ satisfies the corresponding
conditions in the last column of Table \ref{tabu8-geom}.

\end{theorem}

\renewcommand*{\arraystretch}{1.3}
\begin{center}
\begin{table}[h]
    \begin{small}
    \noindent
    \begin{tabular}{|p{1.1cm}|p{5.6cm}|p{4.7cm}|}
              \hline
    class &  normal form & geometric conditions  \\ \hline

  $(U_8)^0$   &  $[U_8]^0_1: [\theta _1 + c_1\theta _2 + c_2\theta _3]_{U_8},\ c_1\ne 0$
    &  $\omega|_V \neq 0$ and $\omega|_{\ell_1+\ell_2} \neq 0$  \\ \cline{2-3}
    &  $[U_8]^0_{\infty}: [\theta _1 + c_2\theta _3]_{U_8}$
    &  $  \omega|_V \neq 0$ and $\omega|_{\ell_1+\ell_2} = 0$ \\ \hline
  $(U_8)^1$ & $[U_8]^1_5: [\pm\theta _2 + c_1\theta _3 + c_2\theta _4]_{U_8}, \ c_2\ne 2c_1$ & $\omega|_V=0$, $\omega|_{\ell_1+\ell_2}\ne 0$ and $L_3=5$    \\ \cline{2-3}

  & $[U_8]^1_{\infty}: [\pm\theta _2 + c_1\theta_3+2c_1\theta_4]_{U_8}$ & $\omega|_V=0$, $\omega|_{\ell_1+\ell_2}\ne 0$ and $L_3=\infty$ \\ \hline
 $(U_8)^{2}$ & $[U_8]^{2}: [\theta _3 + c_1\theta_4+c_2\theta_5]_{U_8}$,\newline $c_1\ne 2$, $c_1\ne -\frac{1}{3}$ &  $\ker\omega = \ell_2$ and $L_3=5$ \\ \hline
 $(U_8)^{3,0}_5$ & $[U_8]^{3,0}_5: [\theta _3- \frac{1}{3}\theta_4+c_1\theta_5+c_2\theta_6]_{U_8}$ & $\ker\omega = \ell_2$ and $L_3=5$ \\ \hline

  $(U_8)^{3,0}_{\infty}$ & $[U_8]^{3,0}_{\infty}: [\theta _3 + 2\theta_4+c_1\theta_5+c_2\theta_7]_{U_8}$ & $\ker\omega = \ell_2$ and $L_3=\infty$ \\ \hline\hline

 & & $\omega\vert_ W = 0$ and $L_{1,2}=\infty$\\ \hline
$(U_8)^{3,1}$  & $[U_8]^{3,1}: [\theta _4 +c\theta_5 ]_{U_8}$
                &  $Lt=L_3=5$  \\ \hline
$(U_8)^4$ & $[U_8]^4: [\pm\theta _5 + c\theta _6]_{U_8}$
                &  $Lt=6$, $L_3=\infty$   \\ \hline
$(U_8)^5$ & $[U_8]^5: [\theta _6 + c\theta _7]_{U_8}$  &  $Lt=7$, $L_3=\infty$   \\ \hline
$(U_8)^6$ & $[U_8]^6: [\pm\theta _7+c\theta_8]_{U_8}$  &  $Lt=8,\;L_3=\infty$   \\ \hline
$(U_8)^7$ & $[U_8]^7: [\theta_8]_{U_8}$     &  $Lt=9,\;L_3=\infty$   \\ \hline
$(U_8)^8$ & $[U_8]^8: [0]_{U_8}$
                &  $N$ is contained in a smooth \newline Lagrangian submanifold   \\ \hline
    \end{tabular}

\medskip

\caption{\small Geometric interpretation of singularity classes of $U_8$. ($W$ is the tangent space to a non-singular 3-dimensional manifold in $(\mathbb R^{2n\geq4}, \omega)$ containing $N\in(U_8)$. The forms $\theta_1,\ldots,\theta_8$ are described in  Theorem ~\ref{u8-baza} on the page \pageref{u8-baza}.)}\label{tabu8-geom}
\end{small}
\end{table}
\end{center}

\begin{xrem}
The idea of the proof of Theorem \ref{geom-cond-u8} is the same as for the proof of Theorem \ref{geom-cond-u7}.
\end{xrem}


\section{Symplectic $U_9$-singularities}\label{sec-u9}

 Denote by $(U_9)$  the class of varieties in  a fixed symplectic space $(\mathbb R^{2n}, \omega )$ which are diffeomorphic to
\begin{equation}
\label{defu9} U_9=\{x\in \mathbb R ^{2n\geq
4}\,:x_1^2+ x_2x_3=x_1x_2+x_3^4=x_{\geq 4}=0\}.\end{equation}

\noindent This is the simple  $1$-dimensional isolated complete intersection singularity $U_9$ (\cite{G}, \cite{AVG}).
Here $N$ is quasi-homogeneous with weights $w(x_1)\!=\!5,\, w(x_2)\!=\!7,\, w(x_3)\!=\!3$.

\medskip

The complete classification of symplectic singularities of $(U_9)$  was obtained using the method of algebraic restrictions.

\begin{theorem}\label{u9-main}
Any submanifold of the symplectic space $(\mathbb R^{2n},\sum_{i=1}^n dp_i \wedge dq_i)$ where $n\geq3$ (respectively $n=2$) which is diffeomorphic to $U_9$ is symplectically equivalent to one and only one of the normal forms $U_9^i, i = 0,1,\cdots ,9$  listed below. The~parameters $c, c_1, c_2, c_3$ of the normal forms are moduli:

\medskip
\begin{small}

\noindent $U_9^0:$ \ $p_1^2 + p_2q_1 = 0, \ \ \pm p_1p_2 +q_1^4= 0, \
\ q_2 = c_1q_1 \mp c_2p_1, \ \ p_{\ge 3} = q_{\ge 3} = 0$;

\smallskip

\noindent $U_9^1:$ \ $p_2^2 \pm p_1q_1 = 0, \ \ p_1p_2 + q_1^4= 0, \
\ q_2 = c_1p_1 + \frac{c_2}{2}q_1^2 \pm \frac{c_3}{3}q_1^3, \ \ p_{\ge 3} = q_{\ge 3} = 0$;

\smallskip

\noindent $U_9^2:$ \ $p_1^2 \pm q_1q_2 = 0, \ \ \pm p_1q_1+q_2^4 = 0, \
\ p_2 = c_1p_1q_2+\frac{c_2}{2}p_1^2, \ \ p_{\ge 3} = q_{\ge 3} = 0$,  $c_1\ne 0$;

\smallskip

\noindent $U_9^{3,0}:$ \ $p_1^2 \pm q_1q_2 = 0, \ \ \pm p_1q_1+q_2^4 = 0, \
\ p_2 = \frac{c_1}{2}p_1^2+c_2p_1q_2^2, \ \ p_{\ge 3} = q_{\ge 3} = 0$,  $c_1\ne 0$;

\smallskip

\noindent $U_9^{4,0}:$ \ $p_1^2 \pm q_1q_2 = 0, \ \ \pm p_1q_1+q_2^4 = 0, \
\ p_2 = c_1p_1q_2^2+\frac{c_2}{2}p_1^2q_2, \ \ p_{\ge 3} = q_{\ge 3} = 0$;

\smallskip

\noindent $U_9^{3,1}:$ \ $p_1^2 + p_2p_3 = 0,  \ p_1p_2 +p_3^4 = 0, \ q_1 =  q_2=0, q_3=-p_1p_3-\frac{c}{2}p_1^2,  \ p_{>3} = q_{>3} = 0$;

\smallskip

\noindent $U_9^{4,1}:$ \ $p_1^2 + p_2p_3 = 0,  \ p_1p_2 +p_3^4 = 0, \ q_1 =  q_2=0, q_3=-\frac{1}{2}p_1^2-c_1p_1p_3^2-c_2p_1p_3^3$, \par  $ \ \ \ \ \  \ p_{>3} = q_{>3} = 0$;

\smallskip

\noindent $U_9^5:$ \ $p_1^2 + p_2p_3 = 0,  \ p_1p_2 +p_3^4 = 0, \ q_1 = q_2=0, q_3=\mp p_1p_3^2-\frac{c}{2}p_1^2p_3,  \ p_{>3} = q_{>3} = 0$;

\smallskip

\noindent $U_9^6:$ \ $p_1^2 + p_2p_3 = 0,  \ p_1p_2 +p_3^4 = 0, \ q_1 = q_2=0, q_3=\mp \frac{1}{2}p_1^2p_3-cp_1p_3^3,  \ p_{>3} = q_{>3} = 0$;

\smallskip

\noindent $U_9^7:$ \ $p_1^2 + p_2p_3 = 0,  \ p_1p_2 +p_3^4 = 0, \  q_1 = q_2=0, q_3 =-p_1p_3^3-\frac{c}{2}p_1^2p_3^2, \ p_{>3} = q_{>3} = 0$;

\smallskip

\noindent $U_9^{8}:$ \ $p_1^2 + p_2p_3 = 0,  \ p_1p_2 +p_3^4 = 0, \  q_1 = q_2=0, q_3 =-\frac{1}{2}p_1^2p_3^2, \ p_{>3} = q_{>3} = 0$;

\smallskip

\noindent $U_9^{9}:$ \ $p_1^2 + p_2p_3 = 0,  \ p_1p_2 +p_3^4 = 0, \
q_{\ge 1} = p_{\ge 4} = 0.$
\end{small}
\end{theorem}


\medskip

\subsection{Distinguishing symplectic classes of $U_9$ by Lagrangian tangency orders}
\label{u9-lagr}
The Lagrangian tangency orders were used to distinguish the symplectic classes of $(U_9)$. A curve $N\in (U_9)$ may be described as a union of two parametrical branches:  $B_1$ and $B_2$. The curve $B_1$ is nonsingular and the curve $B_2$ is singular. Their parametrization  in the coordinate system $(p_1,q_1,p_2,q_2,\cdots,p_n,q_n)$ is presented in the second column of Table \ref{tabu9-lagr}. To characterize the symplectic classes of this singularity we use the following two invariants:
\begin{itemize}
  \item $Lt=Lt(B_1,B_2)=\max\limits _L (\min \{t(B_1,L),t(B_2,L)\}),$
  \item $L_{2}=Lt(B_2)=\max\limits _L\, t(B_2,L).$
\end{itemize}
Here $L$ is a smooth Lagrangian submanifold of the symplectic space.

We can also compare the Lagrangian tangency orders with the respective indices of isotropy.

\begin{theorem}
\label{lagr-u9} A stratified submanifold $N\in (U_9)$ of the symplectic space $(\mathbb R^{2n}, \omega_0)$ with the canonical coordinates $(p_1, q_1, \cdots, p_n, q_n)$ is symplectically equivalent to one and only one of the curves presented in the second column of Table \ref{tabu9-lagr}. The parameters $c, c_1, c_2, c_3$ are moduli. The Lagrangian tangency orders  are presented in  the third and fourth  column of Table \ref{tabu9-lagr}.
\end{theorem}

\renewcommand*{\arraystretch}{1.3}
\begin{center}
\begin{table}[h]

    \begin{small}
    \noindent
    \begin{tabular}{|p{1.1cm}|p{6.7cm} r|c|c|c|c|}
                      \hline
    class &  parametrization of branches  &   & $ind$& $ind_2$ & $Lt$    & $L_2$ \\ \hline

 $(U_9)^0$ & $B_1:(0,0,t,0,0,\cdots )$, &  $c_1\ne 0$ & $0$ & $0$ & $3$  & $5$ \\ \cline{3-7}
     $2n\ge 4$ & $B_2:(\pm t^5,t^3,-t^7,c_1t^3-c_2t^5,0,\cdots )$ &  $c_1=0$ & $0$ & $0$ & $5$  & $5$ \\ \hline

$(U_9)^1$ & $B_1:(t,0,0,c_1t,0,\cdots ),$ &   & $0$ & $0$  & $3$  & $7$ \\
     $2n\ge 4$ & $B_2:(-t^7,\pm t^3,t^5,-c_1 t^7+\frac{c_2}{2}t^6\pm\frac{c_3}{3}t^9,0,\cdots )$ & & & &  &  \\ \hline

$(U_9)^2$ & $B_1:(0,\pm t,0,0,0,\cdots ),$ & $c_1\ne 0$  & $0$ & $0$ & $5$  & $7$ \\
     $2n\ge 4$ & $B_2:(t^5,\mp t^7,c_1t^8+\frac{c_2}{2}t^{10}, t^3,0,\cdots )$ & &  & & &   \\ \hline

$(U_9)^{3,0}$ & $B_1:(0,\pm t,0,0,0,\cdots ),$ & $c_1\ne 0$  & $0$ & $0$ & $5$  & $7$ \\
     $2n\ge 4$ & $B_2:(t^5,\mp t^7,\frac{c_1}{2}t^{10}+c_2t^{11}, t^3,0,\cdots )$ & &  & & &   \\ \hline

$(U_9)^{4,0}$ & $B_1:(0,\pm t,0,0,0,\cdots ),$ &   & $0$ & $0$  & $5$  & $7$ \\
     $2n\ge 4$ & $B_2:(t^5,\mp t^7,c_1t^{11}+\frac{c_2}{2}t^{13}, t^3,0,\cdots )$ & &  & & &   \\ \hline

$(U_9)^{3,1}$ & $B_1:(0,0,t,0,0,0,\cdots ),$ &   & $1$ & $1$ & $8$  & $8$ \\
     $2n\ge 6$ & $B_2:(t^5,0,-t^7,0, t^3,-t^8-\frac{c}{2}t^{10},0,\cdots )$ & & & &  &  \\ \hline

$(U_9)^{4,1}$ & $B_1:(0,0,t,0,0,0,\cdots ),$ &   & $1$ & $1$ & $10$  & $10$ \\
     $2n\ge 6$ & $B_2:(t^5,0,-t^7,0, t^3,-\frac{1}{2}t^{10}-c_1t^{11}-c_2t^{14},0,\cdots )$ & & & &  &  \\ \hline

$(U_9)^5$ & $B_1:(0,0,t,0,0,0,\cdots ),$ &   & $2$ & $2$ & $11$  & $11$ \\
     $2n\ge 6$ & $B_2:(t^5,0,-t^7,0, t^3,\mp t^{11}-\frac{c}{2}t^{13},0,\cdots )$ & & & &  &  \\ \hline

$(U_9)^6$ & $B_1:(0,0,t,0,0,0,\cdots ),$ &   & $2$ & $2$ & $13$  & $13$ \\
     $2n\ge 6$ & $B_2:(t^5,0,-t^7,0, t^3,\mp \frac{1}{2}t^{13}-ct^{14},0,\cdots )$ & & & &  &  \\ \hline

$(U_9)^7$ & $B_1:(0,0,t,0,0,0,\cdots ),$ &   & $3$ & $\infty$ & $14$  & $\infty$ \\
     $2n\ge 6$ & $B_2:(t^5,0,-t^7,0, t^3,-t^{14}-\frac{c}{2}t^{16},0,\cdots )$ & & & &  &  \\ \hline

$(U_9)^{8}$ & $B_1:(0,0,t,0,0,0,\cdots ),$ &   & $3$ & $\infty$ & $16$  & $\infty$ \\
     $2n\ge 6$ & $B_2:(t^5,0,-t^7,0, t^3,-\frac{1}{2}t^{16},0,\cdots )$ &  & & & &  \\ \hline

$(U_9)^{9}$ & $B_1:(0,0,t,0,0,0,\cdots ),$ &   & $\infty$ & $\infty$  & $\infty$ & $\infty$ \\
     $2n\ge 6$ & $B_2:(t^5,0,-t^7,0, t^3,0,0,\cdots )$ & &  & &  &  \\ \hline

\end{tabular}

\medskip

\caption{\small The Lagrangian tangency orders for symplectic classes of the $U_9$ singularity.}\label{tabu9-lagr}
\end{small}
\end{table}
\end{center}
\medskip

\subsection{Geometric conditions for the classes $(U_9)^i$}
\label{u9-geom_cond} $ \ $



Let $N\in (U_9)$. Denote by $W$ the tangent space at $0$ to some (and then any) non-singular $3$-manifold containing $N$. We can define the following subspaces of this space:

  $\ell_1$ -- the tangent line at $0$ to the nonsingular branch $B_1$,

  $\ell_2$ -- the tangent line at $0$ to the singular branch $B_2$,

  $V$ -- the $2$-space tangent at $0$ to the singular branch $B_2$.

    For $N=U_9=$(\ref{defu9}) it is easy to calculate that
  $W\!=\!\Span (\partial /\partial x_1,\partial /\partial x_2,\partial /\partial x_3)$,  \newline and
$\ell\!_1=\!\Span (\partial /\partial x_2), \ \ell\!_2=\!\Span (\partial /\partial x_3), \
V\!=\!\Span (\partial /\partial x_1,\partial /\partial x_3)$.

The classes $(U_9)^i$ satisfy special conditions in terms of the restriction $\omega\vert_ W $, where $\omega $ is the symplectic form.

\begin{theorem}
\label{geom-cond-u9} For any stratified submanifold $N\in (U_9)$ of the symplectic space $(\mathbb R^{2n}, \omega )$ belonging to the class $(U_9)^i$  the couple $(N, \omega )$ satisfies the corresponding conditions in the last column of Table \ref{tabu9-geom}.
\end{theorem}


\renewcommand*{\arraystretch}{1.3}
\setlength{\tabcolsep}{1.2mm}
\begin{center}
\begin{table}[h]
    \begin{small}
    \noindent
    \begin{tabular}{|p{1.0cm}|p{5.5cm}|p{5.3cm}|}
            \hline
    class &  normal form & geometric conditions  \\ \hline

  $(U_9)^0$   &  $[U_9]^0_0: [\pm\theta _1 + c_1\theta _2 + c_2\theta _3]_{U_9},\ c_1\ne 0$
    &  $\omega|_V \neq 0$ and $\omega|_{\ell_1+\ell_2} \neq 0$  \\ \cline{2-3}
    &  $[U_9]^0_1: [\pm\theta _1 + c_2\theta _3]_{U_9}$
    &  $  \omega|_V \neq 0$ and $\omega|_{\ell_1+\ell_2} = 0$ \\ \hline
  $(U_9)^1$ & $[U_9]^1: [\pm\theta _2 + c_1\theta _3 + c_2\theta _4+ c_3\theta _6]_{U_9}$ & $\omega|_V=0$, $\omega|_{\ell_1+\ell_2} = 0$ and $\ker\omega \neq \ell_2$    \\ \hline

 $(U_9)^{2}$ & $[U_9]^2: [\pm\theta _3 + c_1\theta_4+c_2\theta_5]_{U_9},\ c_1\ne 0$ &  \\ \cline{1-2}
 $(U_9)^{3,0}$ & $[U_9]^{3,0}: [\pm\theta _3 + c_1\theta_5+c_2\theta_6]_{U_9},\ c_1\ne 0$ & $\omega|_V=0$ and $\ker\omega = \ell_2$ \\ \cline{1-2}
 $(U_9)^{4,0}$ & $[U_9]^{4,0}: [\pm\theta _3 + c_1\theta_6+c_2\theta_7]_{U_9}$ &  \\ \hline\hline

 & & $\omega\vert_ W = 0$\\ \hline
$(U_9)^{3,1}$  & $[U_9]^{3,1}: [\theta _4 +c\theta_5 ]_{U_9}$
                &  $Lt=L_2=8$  \\ \hline
$(U_9)^{4,1}$ & $[U_9]^{4,1}: [\theta _5 + c_1\theta _6+c_2\theta _8]_{U_9}$
                &  $Lt=L_2=10$   \\ \hline
$(U_9)^5$ & $[U_9]^5: [\pm\theta _6 + c\theta _7]_{U_9}$
                &  $Lt=L_2=11$   \\ \hline
$(U_9)^6$ & $[U_9]^6: [\pm\theta _7 + c\theta _8]_{U_9}$  &  $Lt=L_2=13$   \\ \hline
$(U_9)^7$ & $[U_9]^7: [\theta _8+c\theta_9]_{U_9}$  &  $Lt=14,\;L_2=\infty$   \\ \hline
$(U_9)^8$ & $[U_9]^8: [\theta_9]_{U_9}$     &  $Lt=16,\;L_2=\infty$   \\ \hline
$(U_9)^9$ & $[U_9]^9: [0]_{U_9}$
                &   $N$ is contained in a smooth \newline Lagrangian submanifold   \\ \hline
    \end{tabular}

\medskip

\caption{\small Geometric characterization of symplectic classes of the $U_9$ singularity. (The forms $\theta_1,\ldots,\theta_9$ are described in Theorem~\ref{u9-baza} on the page \pageref{u9-baza}.)}\label{tabu9-geom}

\end{small}
\end{table}
\end{center}
\begin{xrem}
The idea of the proof of Theorem \ref{geom-cond-u9} is the same as for the proof of Theorem \ref{geom-cond-u7}.
\end{xrem}

\medskip

\section{Proofs}\label{proofs}

\subsection{The method of algebraic restrictions}
\label{method}

In this section we present only basic notions and facts on the method of algebraic restrictions, which is a very powerful tool for the symplectic classification. The details of the method  can be found in \cite{DJZ2}.

  Given a germ of a non-singular manifold $M$ denote by $\Lambda ^p(M)$ the space of all germs at $0$ of differential $p$-forms on $M$. Given a subset $N\subset M$ introduce the following subspaces of $\Lambda ^p(M)$:
$$\Lambda ^p_N(M) = \{\omega \in \Lambda ^p(M): \ \ \omega (x)=0 \ \text {for any} \ x\in N \};$$
$$\mathcal A^p_0(N, M) = \{\alpha  + d\beta : \ \ \alpha \in \Lambda _N^p(M), \ \beta \in \Lambda _N^{p-1}(M)\}.$$

\smallskip

\begin{defn}
\label{main-def} Let $N$ be the germ of a subset of $M$ and let
$\omega \in \Lambda ^p(M)$. The \emph{algebraic restriction} of
$\omega $ to $N$ is the equivalence class of $\omega $ in $\Lambda
^p(M)$, where the equivalence is as follows: $\omega $ is
equivalent to $\widetilde \omega $ if $\omega - \widetilde \omega
\in \mathcal A^p_0(N, M)$.
\end{defn}

\begin{nt} The algebraic restriction of the germ of
a $p$-form $\omega $ on $M$ to the germ of a subset $N\subset M$
will be denoted by $[\omega ]_N$. By writing $[\omega ]_N=0$ (or
saying that $\omega $ has zero algebraic restriction to $N$) we
mean that $[\omega ]_N = [0]_N$, i.e. $\omega \in A^p_0(N, M)$.
\end{nt}
\medskip

\begin{defn}Two algebraic restrictions
$[\omega ]_N$ and $[\widetilde \omega ]_{\widetilde N}$ are called \emph{diffeomorphic} if there exists the germ of a diffeomorphism $\Phi:
\widetilde M\to M$ such that $\Phi(\widetilde N)=N$ and  $\Phi ^*([\omega ]_N) =[\widetilde \omega ]_{\widetilde N}$.
\end{defn}

\smallskip

The method of algebraic restrictions applied to singular
quasi-homogeneous subsets is based on the following theorem.

\begin{theorem}[Theorem A in \cite{DJZ2}] \label{thm A}
Let $N$ be the germ of a quasi-homogeneous subset of $\mathbb
R^{2n}$. Let $\omega _0, \omega _1$ be germs of symplectic forms
on $\mathbb R^{2n}$ with the same algebraic restriction to $N$.
There exists a local diffeomorphism $\Phi $ such that $\Phi (x) =
x$ for any $x\in N$ and $\Phi ^*\omega _1 = \omega _0$.

Two germs of quasi-homogeneous subsets $N_1, N_2$ of a fixed
symplectic space $(\mathbb R^{2n}, \omega )$ are symplectically
equivalent if and only if the algebraic restrictions of the
symplectic form $\omega $ to $N_1$ and $N_2$ are diffeomorphic.

\end{theorem}

\medskip

Theorem \ref{thm A} reduces the problem of symplectic
classification of germs of singular quasi-homogeneous subsets to
the problem of diffeomorphic classification of algebraic
restrictions of the germ of the symplectic form to the germs of
singular quasi-homogeneous subsets.

\smallskip

The geometric meaning of the zero algebraic restriction is explained
by the following theorem.

\begin{theorem}[Theorem {\bf B} in \cite{DJZ2}] \label{thm B}  {\it The germ of a quasi-homogeneous
 set $N$  of a symplectic space
$(\mathbb R^{2n}, \omega )$ is contained in a non-singular
Lagrangian submanifold if and only if the symplectic form $\omega$ has zero algebraic restriction to $N$.}\
\end{theorem}

In the remainder of this paper we use the following notations:

\smallskip

\noindent $\bullet$  $\algrestall $: \ \ the vector space
consisting of  the algebraic restrictions of germs of all $2$-forms on
$\mathbb R^{2n}$ to the germ of a subset $N\subset \mathbb
R^{2n}$;

\smallskip

\noindent $\bullet$  $\algrestclosed$: \ \ the subspace of
$\algrestall $ consisting of the algebraic restrictions of germs of
all closed $2$-forms on $\mathbb R^{2n}$ to $N$;

\smallskip

\noindent $\bullet$ $\algrest $: \ \ the open set in
$\algrestclosed$ consisting of the algebraic restrictions of germs of
all symplectic $2$-forms on $\mathbb R^{2n}$ to $N$.

\medskip

To obtain a classification of the algebraic restrictions we use the following proposition.
\begin{proposition}\label{elimin1}
Let $a_1, \cdots, a_p$ be a quasi-homogeneous basis of
quasi-degrees
$\delta_1\le\cdots\le\delta_s<\delta_{s+1}\le\cdots\le\delta_p$ of
the space of algebraic restrictions of closed $2$-forms to quasi-homogeneous subset $N$.
Let $a=\sum_{j=s}^p c_ja_j$, where $c_j\in \mathbb R$ for
$j=s,\cdots, p$ and $c_s\ne 0$.

If there exists a tangent quasi-homogeneous vector field $X$ over
$N$ such that $\mathcal L_Xa_s=ra_k$ for $k>s$ and $r\ne 0$ then
$a$ is diffeomorphic to $\sum_{j=s}^{k-1} c_ja_j+\sum_{j=k+1}^p
b_j a_j$, for some $b_j\in \mathbb R, \ j=k+1,\cdots,p$.
\end{proposition}

Proposition \ref{elimin1} is a modification of  Theorem 6.13 formulated and proved in  \cite{D}. It was formulated for algebraic restrictions to a parameterized curve but we can generalize this theorem for any quasi-homogeneous subset $N$. The proofs of the cited theorem and Proposition \ref{elimin1}  are based on the Moser homotopy method.

\medskip

For calculating discrete invariants we use the following propositions.

\begin{proposition}[\cite{DJZ2}]\label{sm}
The symplectic multiplicity  of the germ of a quasi-homo\-ge\-neous subset $N$ in a symplectic space is equal to the codimension of the orbit of the algebraic restriction $[\omega ]_N$ with respect to the group of local diffeomorphisms preserving $N$  in the space of algebraic restrictions of closed  $2$-forms to $N$.
\end{proposition}

\begin{proposition}[\cite{DJZ2}]\label{ii}
The index of isotropy  of the germ of a quasi-homogeneous subset $N$ in a symplectic space $(\mathbb R^{2n}, \omega )$ is equal to the maximal order of vanishing of closed $2$-forms representing the algebraic restriction $[\omega ]_N$.
\end{proposition}

\begin{proposition}[\cite{D}]\label{lto}
Let $f$ be the germ of a quasi-homogeneous curve such that the algebraic restriction of a symplectic form to it can be represented by a closed $2$-form vanishing at $0$. Then the Lagrangian tangency order of the germ of a quasi-homogeneous curve $f$ is the maximum of the order of vanishing on $f$ over all $1$-forms $\alpha$ such that $[\omega]_f=[d\alpha]_f$
\end{proposition}

\subsection{Proofs for  $U_7$ singularity}
\subsubsection{Algebraic restrictions to  $U_7$ and their classification}\label{u7-class}

One has the following relations for $(U_7)$-singularities:
\begin{equation}
[x_1^2+x_2x_3]_{U_7}=0.
\label{u01}
\end{equation}
\begin{equation}
[x_1x_2+x_3^3]_{U_7}=0,
\label{u02}
\end{equation}
\begin{equation}
[d(x_1^2+x_2x_3)]_{U_7}=[2x_1dx_1+ x_2dx_3+x_3 dx_2]_{U_7}=0
\label{u71}
\end{equation}
\begin{equation}
[d(x_1x_2+x_3^3)]_{U_7}=[x_1dx_2+x_2dx_1+3x_3^2dx_3]_{U_7}=0
\label{u72}
\end{equation}
Multiplying these relations by suitable $1$-forms and $2$-forms we obtain the relations towards calculating $[\Lambda^2(\mathbb R^{2n})]_N$ for $N=U_7$. 

\begin{proposition}
\label{u7-all}
The space $[\Lambda ^{2}(\mathbb R^{2n})]_{U_7}$ is an $8$-dimensional vector space spanned by the algebraic restrictions to $U_7$ of the $2$-forms

$\theta _1= dx_1\wedge dx_3, \;\; \theta _2=dx_2\wedge dx_3,\;\; \theta_3 = dx_1\wedge dx_2,$


  $\theta _4 = x_3dx_1\wedge dx_3,\;\; \theta _5 = x_1dx_1\wedge dx_3,$\;\; $\sigma = x_1 dx_2\wedge dx_3,$


$\theta _6= x_3^2 dx_1\wedge dx_3$,\;\; $\theta _7= x_1x_3 dx_1\wedge dx_3$.
\end{proposition}

Proposition \ref{u7-all} and results of Section \ref{method}  imply the following description of the space $[Z ^2(\mathbb R^{2n})]_{U_7}$ and the manifold $[{\rm Symp} (\mathbb R^{2n})]_{U_7}$.

\begin{theorem} \label{u7-baza}
 The space $[Z^2(\mathbb R^{2n})]_{U_7}$ is a $7$-dimensional vector space
 spanned by the algebraic restrictions to $U_7$  of the quasi-homogeneous $2$-forms $\theta_i$  of degree $\delta_i$

\smallskip

$\theta _1= dx_1\wedge dx_3,\;\;\;\delta_1=7,$


$\theta _2=dx_2\wedge dx_3,\;\;\;\delta_2=8,$


$\theta_3 = dx_1\wedge dx_2,\;\;\;\delta_3=9,$


  $\theta _4 = x_3dx_1\wedge dx_3,\;\;\;\delta_4=10,$


  $\theta _5 = x_1dx_1\wedge dx_3,\;\;\;\delta_5=11,$


$\theta _6= x_3^2 dx_1\wedge dx_3,\;\;\;\delta_6=13, $


$\theta _7= x_1x_3 dx_1\wedge dx_3,\;\;\;\delta_7=14$.
\smallskip

If $n\ge 3$ then $[{\rm Symp} (\mathbb R^{2n})]_{U_7} = [Z^2(\mathbb R^{2n})]_{U_7}$. The manifold $[{\rm Symp} (\mathbb R^{4})]_{U_7}$ is an open part of the $7$-space $[Z^2 (\mathbb R^{4})]_{U_7}$ consisting of algebraic restrictions of the form $[c_1\theta _1 + \cdots + c_7\theta _7]_{U_7}$ such that $(c_1,c_2,c_3)\ne (0,0,0)$.
\end{theorem}

\begin{theorem}
\label{klasu7} $ \ $

\smallskip

\noindent (i) \ Any algebraic restriction in $[Z ^2(\mathbb R^{2n})]_{U_7}$ can be brought by a symmetry of $U_7$ to one of the normal forms $[U_7]^i$ given in the second column of Table \ref{tabu7}.

\smallskip

\noindent (ii)  \ The codimension in $[Z ^2(\mathbb R^{2n})]_{U_7}$ of the singularity class corresponding to the normal form $[U_7]^i$ is equal to $i$, the symplectic multiplicity and the index of isotropy are given in the fourth and fifth columns of Table  \ref{tabu7}.

\smallskip

\noindent (iii) \ The singularity classes corresponding to the normal forms are disjoint.

\smallskip

\noindent (iv) \ The parameters $c, c_1, c_2$ of the normal forms $[U_7]^i$ are moduli.

\end{theorem}

\renewcommand*{\arraystretch}{1.3}
\begin{center}
\begin{table}[h]
    \begin{small}
    \noindent
    \begin{tabular}{|p{3cm}|p{5.8cm}|c|c|c|}
                      \hline
    symplectic class &   normal forms for algebraic restrictions    & cod & $\mu ^{\rm sym}$ &  ind  \\ \hline
   $(U_7)^0$ \;\;  $(2n\ge 4)$ & $[U_7]^0: [\theta _1 + c_1\theta _2 + c_2\theta _3]_{U_7}$,\;
                          &  $0$ & $2$ & $0$  \\  \hline
  $(U_7)^1$ \;\; $(2n\ge 4)$ & $[U_7]^1: [\pm\theta _2 + c_1\theta _3 + c_2\theta _4]_{U_7}$
                                         &  $1$ & $3$ & $0$ \\ \hline
    $(U_7)^2$ \;\; $(2n\ge 4)$& $[U_7]^2: [\theta _3 + c_1\theta_4+c_2\theta_5]_{U_7}$,\;
              & $2$ & $4$ & $0$     \\ \hline
    $(U_7)^3$ \;\; $(2n\ge 6)$ & $[U_7]^3: [\pm\theta _4 + c\theta _5]_{U_7}$ &
                $3$ & $4$ & $1$  \\ \hline
    $(U_7)^4$ \;\; $(2n\ge 6)$ & $[U_7]^4: [\theta _5 + c\theta _6]_{U_7}$
                                         &  $4$ & $5$ & $1$    \\ \hline
    $(U_7)^5$ \;\; $(2n\ge 6)$ & $[U_7]^5: [\theta _6 + c\theta _7]_{U_7}$
                                         &  $5$ & $6$ & $2$ \\ \hline
     $(U_7)^6$ \;\; $(2n\ge 6)$ & $[U_7]^6: [\pm\theta _7]_{U_7}$
                                         &  $6$ & $6$ & $2$   \\ \hline
    $(U_7)^7$ \;\; $(2n\ge 6)$ & $[U_7]^7: [0]_{U_7}$ &  $7$ & $7$ & $\infty $ \\ \hline
\end{tabular}

\smallskip

\caption{\small Classification of symplectic $U_7$ singularities.\newline
$cod$ -- codimension of the classes;  $\mu ^{sym}$-- symplectic multiplicity; \newline  $ind$ -- the index of isotropy.}\label{tabu7}
\end{small}
\end{table}
\end{center}


\noindent In the first column of Table \ref{tabu7}   we denote by $(U_7)^i$ a subclass of $(U_7)$ consisting of $N\in (U_7)$ such that the algebraic restriction $[\omega ]_N$ is diffeomorphic to some algebraic restriction of the normal form $[U_7]^i$, where $i$ is the codimension of the class.

The proof of Theorem \ref{klasu7} is presented in Section \ref{u7-proof}.

\subsubsection{Symplectic normal forms}
\label{u7-normal}

Let us transfer the normal forms $[U_7]^i$  to symplectic normal forms
. Fix a family $\omega ^i$ of symplectic forms on $\mathbb R^{2n}$ realizing the family $[U_7]^i$ of algebraic restrictions.  We can fix, for example,

\smallskip
\begin{small}

\noindent $\omega ^0 = \theta _1 + c_1\theta _2 + c_2\theta _3 +
dx_2\wedge dx_4 + dx_5\wedge dx_6 + \cdots + dx_{2n-1}\wedge
dx_{2n};$

\smallskip

\noindent $\omega ^1 = \pm \theta _2 + c_1\theta _3 + c_2\theta _4 +
dx_1\wedge dx_4 + dx_5\wedge dx_6 + \cdots + dx_{2n-1}\wedge
dx_{2n}; $

\smallskip

\noindent $\omega ^2 = \theta_3+ c_1 \theta _4 + c_2 \theta _5 + dx_3\wedge dx_4 +
dx_5\wedge dx_6 + \cdots + dx_{2n-1}\wedge dx_{2n};$

\smallskip

\noindent $\omega ^3 = \pm\theta _4 + c\theta _5 + dx_1\wedge dx_4 +
dx_2\wedge dx_5 + dx_3\wedge dx_6 + dx_7\wedge dx_8 + \cdots +
dx_{2n-1}\wedge dx_{2n};$

\smallskip

\noindent $\omega ^4 = \theta _5 + c\theta _6 + dx_1\wedge dx_4 +
dx_2\wedge dx_5 + dx_3\wedge dx_6 + dx_7\wedge dx_8 +\cdots + dx_{2n-1}\wedge dx_{2n};$

\smallskip

\noindent $\omega ^5 = \theta _6 + c\theta _7+ dx_1\wedge dx_4 + dx_2\wedge dx_5 +
dx_3\wedge dx_6 + dx_7\wedge dx_8 + \cdots + dx_{2n-1}\wedge
dx_{2n};$

\smallskip

\noindent $\omega ^6 = \pm\theta _7+ dx_1\wedge dx_4 + dx_2\wedge dx_5 +
dx_3\wedge dx_6 + dx_7\wedge dx_8 + \cdots + dx_{2n-1}\wedge
dx_{2n};$

\smallskip

\noindent $\omega ^7 = dx_1\wedge dx_4 + dx_2\wedge dx_5 +
dx_3\wedge dx_6 + dx_7\wedge dx_8 + \cdots + dx_{2n-1}\wedge
dx_{2n}.$
\end{small}

\medskip
 Let $\omega_0 = \sum_{i=1}^m dp_i \wedge dq_i$, where $(p_1,q_1,\cdots,p_n,q_n)$ is the coordinate system on $\mathbb R^{2n}, n\ge 3$ (resp. $n=2$). Fix, for $i=0,1,\cdots ,7$ (resp. for $i = 0,1,2)$ a family $\Phi ^i$ of local diffeomorphisms which bring the family of symplectic forms $\omega ^i$ to the symplectic form $\omega_0$: $(\Phi ^i)^*\omega ^i = \omega_0$. Consider the families $U_7^i = (\Phi ^i)^{-1}(U_7)$. Any stratified submanifold of the symplectic space $(\mathbb R^{2n},\omega_0)$ which is diffeomorphic to $U_7$ is symplectically equivalent to one and only one of the normal forms $U_7^i, i = 0,1,\cdots ,7$ (resp. $i= 0,1,2$) presented in Theorem \ref{u7-main}. By Theorem \ref{klasu7} we obtain that  parameters $c,c_1,c_2$ of the normal forms are moduli.

 \medskip

\subsubsection{Proof of Theorem \ref{klasu7}}
\label{u7-proof}

In our proof we use vector fields tangent to $N\in U_7$. Any vector fields tangent to $N\in U_7$ can be described as $V=g_1E+g_2\mathcal{H}$ where $E$ is the Euler vector field and $\mathcal{H}$ is a Hamiltonian vector field and $g_1, g_2$ are functions. It was shown in \cite{DT1} (Prop. 6.13) that the action of a Hamiltonian vector field on the algebraic restriction  of a closed 2-form to any 1-dimensional complete intersection is trivial.

\medskip

\noindent The germ of a vector field tangent to $U_7$ of non trivial action on algebraic restrictions of closed 2-forms to  $U_7$ may be described as a linear combination of germs of vector fields:   $X_0=E,\, X_1=x_3E,\, X_2=x_1E,\, X_3=x_2E,\, X_4=x_3E^2,\, X_5=x_1x_3E$, where $E$ is the Euler vector field    \begin{equation}
\label{eu7}
E=4 x_1 \partial /\partial x_1+5 x_2 \partial /\partial x_2+3 x_3 \partial /\partial x_3.
\end{equation}

\begin{proposition} \label{u7-infinitesimal}

The infinitesimal action of germs of quasi-homogeneous vector
fields tangent to $N\in (U_7)$ on the basis of the vector space of
algebraic restrictions of closed $2$-forms to $N$ is presented in
Table \ref{infini-u7}.

\setlength{\tabcolsep}{2.4mm}
\renewcommand*{\arraystretch}{1.3}
\begin{small}
\begin{table}[h]
\begin{center}
\begin{tabular}{|l|r|r|r|r|r|r|r|}

 \hline

  $\mathcal L_{X_i} [\theta_j]$ & $[\theta_1]$   &   $[\theta_2]$ &   $[\theta_3]$ & $[\theta_4]$   & $[\theta_5]$   & $[\theta_6]$   & $[\theta_7]$  \\ \hline

  $X_0=E$ & $7 [\theta_1]$ & $8 [\theta_2]$ & $9 [\theta_3]$ & $10 [\theta_4]$ & $11 [\theta_5]$ & $13 [\theta_6]$ & $14 [\theta_7]$  \\ \hline

  $X_1=x_3E$ & $10[\theta_4]$ & $-22 [\theta_5]$ & $ [0]$ & $13 [\theta_6]$ & $14[\theta_7]$  & $[0]$ & $[0]$   \\  \hline

  $X_2=x_1E$ & $11[\theta_5]$  & $ [0]$ & $-39 [\theta_6]$ & $14[\theta_7]$ & $ [0]$ & $[0]$ & $[0]$  \\  \hline

  $X_3=x_2 E$ & $[0]$  & $-78[\theta_6]$ & $-84[\theta_7]$  & $[0]$ & $[0]$ & $[0]$  & $[0]$   \\  \hline

  $X_4=x_3^2E$ & $13 [\theta_6]$ & $-28[\theta_7]$  & $[0]$ & $[0]$ & $[0]$ & $[0]$ & $[0]$  \\ \hline

  $X_5=x_1x_3 E$ & $14[\theta_7]$  & $[0]$ & $[0]$  & $[0]$ & $[0]$   & $[0]$ & $[0]$ \\  \hline
\end{tabular}
\end{center}

\caption{\small Infinitesimal actions on algebraic restrictions of closed \newline
2-forms to   $U_7$. ($E$ is defined as in (\ref{eu7}).)}\label{infini-u7}
\end{table}
\end{small}
\end{proposition}

Let $\mathcal{A}=[c_1 \theta_1+c_2 \theta_2+c_3 \theta_3+c_4 \theta_4+c_5 \theta_5+c_6 \theta_6 +c_7 \theta_7]_{U_7}$
be the algebraic restriction of a symplectic form $\omega$.

\smallskip

The first statement of Theorem \ref{klasu7} follows from the following lemmas.

\begin{lemma}
\label{u7lem0} If \;$c_1\ne 0$\; then the algebraic restriction
$\mathcal{A}=[\sum_{k=1}^7 c_k \theta_k]_{U_7}$
can be reduced by a symmetry of $U_7$ to an algebraic restriction $[\theta_1+\widetilde{c}_2 \theta_2+\widetilde{c}_3 \theta_3]_{U_7}$.
\end{lemma}




\begin{lemma}
\label{u7lem1} If \;$c_1\!=\!0$\;and $c_2\ne\!0$ then the algebraic restriction  $\mathcal{A}$ 
can be reduced by a symmetry of $U_7$ to an algebraic restriction $[\pm\theta_2+\widetilde{c}_3 \theta_3+\widetilde{c}_4 \theta_4]_{U_7}$.
\end{lemma}

\begin{lemma}
\label{w7lem2} If $c_1\!=\!c_2\!=\!0$ and $c_3\!\neq\!0$ then the algebraic restriction $\mathcal{A}$
can be reduced by a symmetry of \, $U_7$ to an algebraic restriction $[\theta_3+\widetilde{c}_4 \theta_4+\widetilde{c}_5 \theta_5]_{U_7}$.
\end{lemma}

\begin{lemma}
\label{u7lem3} If \;$c_1=c_2=c_3=0$ and $c_4\ne 0$\; then the algebraic restriction $\mathcal{A}$
can be reduced by a symmetry of \, $U_7$ to an algebraic restriction   $[\pm\theta_4+\widetilde{c}_5 \theta_5]_{U_7}$.
\end{lemma}

\begin{lemma}
\label{u7lem4} If \;$c_1=0,\ldots,c_4=0$ and $c_5\ne 0$,\; then the algebraic restriction $\mathcal{A}$
can be reduced by a symmetry of \, $U_7$ to an algebraic restriction $[\theta_5+\widetilde{c}_6 \theta_6]_{U_7}$.
\end{lemma}

\begin{lemma}
\label{u7lem5} If \;$c_1=0,\ldots,c_5=0$ and $c_6\ne 0$\; then the algebraic restriction $\mathcal{A}$
can be reduced by a symmetry of \, $U_7$ to an algebraic restriction $[\theta_6+\widetilde{c}_7 \theta_7]_{U_7}$.
\end{lemma}

\begin{lemma}
\label{u7lem6} If \;$c_1=0,\ldots,c_6=0$ and $c_7\ne 0$\; then the algebraic restriction $\mathcal{A}$
can be reduced by a symmetry of \, $U_7$ to an algebraic restriction $[\pm\theta_7]_{U_7}$.
\end{lemma}


The proofs of Lemmas \ref{u7lem0} -- \ref{u7lem6} are similar and based on Table \ref{infini-u7} and Proposition \ref{elimin1}.

\medskip

Statement $(ii)$ of Theorem \ref{klasu7} follows from the conditions
in the proof of part $(i)$ (the codimension) and  from Theorem \ref{thm B} and Proposition \ref{sm} (the symplectic multiplicity) and Proposition \ref{ii} (the index of isotropy).

\medskip

To prove statement $(iii)$ of Theorem \ref{klasu7} we have to show that singularity classes corresponding to normal forms are disjoint. It is enough to notice that  the singularity classes  can be distinguished by geometric conditions. 

\medskip

To prove statement $(iv)$ of Theorem \ref{klasu7} we have to show that the parameters $c, c_1, c_2$ are moduli in the normal forms. The proofs are very similar in all cases. We consider as an example
the normal form with two parameters $[\theta_1+c_1\theta_2+c_2\theta_3]_{U_7}$. From Table \ref{infini-u7} we see that the tangent space to the orbit
of $[\theta_1+c_1\theta_2+c_2\theta_3]_{U_7}$ at $[\theta_1+c_1\theta_2+c_2\theta_3]_{U_7}$ is spanned by the linearly independent algebraic restrictions
$[7\theta_1+8c_1\theta_2+9c_2\theta_3]_{U_7}$, $[\theta_4]_{U_7},[\theta_5]_{U_7}, [\theta_6]_{U_7}$ and $[\theta_7]_{U_7}.$ Hence, the algebraic restrictions
$[\theta_2]_{U_7}$ and $[\theta_3]_{U_7}$ do not belong to it. Therefore, the parameters $c_1$ and $c_2$ are independent moduli in the normal form
$[\theta_1+c_1\theta_2+c_2\theta_3]_{U_7}$.

\bigskip

\subsection{Proofs for  $U_8$ singularity}
\subsubsection{Algebraic restrictions to  $U_8$ and their classification}\label{u8-class}

One has the following relations for $(U_8)$-singularities:
\begin{equation}
[x_1^2+x_2x_3]_{U_8}=0.
\label{u801}
\end{equation}
\begin{equation}
[x_1x_2+x_1x_3^2]_{U_8}=0,
\label{u802}
\end{equation}
\begin{equation}
[d(x_1^2+x_2x_3)]_{U_8}=[2x_1dx_1+ x_2dx_3+x_3 dx_2]_{U_8}=0
\label{u81}
\end{equation}
\begin{equation}
[d(x_1x_2+x_1x_3^2)]_{U_8}=[x_1dx_2+x_2dx_1+x_3^2dx_1+2x_1x_3dx_3]_{U_8}=0
\label{u82}
\end{equation}
Multiplying these relations by suitable $1$-forms and $2$-forms we obtain the relations towards calculating $[\Lambda^2(\mathbb R^{2n})]_N$ for $N=U_8$. 

\begin{proposition}
\label{u8-all}
The space $[\Lambda ^{2}(\mathbb R^{2n})]_{U_8}$ is a $9$-dimensional vector space spanned by the algebraic restrictions to $U_8$ of the $2$-forms

$\theta _1= dx_1\wedge dx_3, \;\; \theta _2=dx_2\wedge dx_3,\;\; \theta_3 = dx_1\wedge dx_2,$


  $\theta _4 = x_3dx_1\wedge dx_3,\;\; \theta _5 = x_1dx_1\wedge dx_3,$\;\; $\theta _6= x_3^2 dx_1\wedge dx_3$, \; \; $\sigma = x_1 dx_2\wedge dx_3,$


 $\theta _7= x_1x_3 dx_1\wedge dx_3$,\;\; $\theta _8= x_3^3 dx_1\wedge dx_3$.
\end{proposition}

Proposition \ref{u8-all} and results of Section \ref{method}  imply the following description of the space $[Z ^2(\mathbb R^{2n})]_{U_8}$ and the manifold $[{\rm Symp} (\mathbb R^{2n})]_{U_8}$.

\begin{theorem} \label{u8-baza}
 The space $[Z^2(\mathbb R^{2n})]_{U_8}$ is an $8$-dimensional vector space
 spanned by the algebraic restrictions to $U_8$  of the quasi-homogeneous $2$-forms $\theta_i$  of degree $\delta_i$

\smallskip

$\theta _1= dx_1\wedge dx_3,\;\;\;\delta_1=5,$


$\theta _2=dx_2\wedge dx_3,\;\;\;\delta_2=6,$


$\theta_3 = dx_1\wedge dx_2,\;\;\;\delta_3=7,$


  $\theta _4 = x_3dx_1\wedge dx_3,\;\;\;\delta_4=7,$


  $\theta _5 = x_1dx_1\wedge dx_3,\;\;\;\delta_5=8,$


$\theta _6= x_3^2 dx_1\wedge dx_3,\;\;\;\delta_6=9, $


$\theta _7= x_1x_3 dx_1\wedge dx_3,\;\;\;\delta_7=10$,

$\theta _8= x_3^3 dx_1\wedge dx_3,\;\;\;\delta_8=11$.

\smallskip

If $n\ge 3$ then $[{\rm Symp} (\mathbb R^{2n})]_{U_8} = [Z^2(\mathbb R^{2n})]_{U_8}$. The manifold $[{\rm Symp} (\mathbb R^{4})]_{U_8}$ is an open part of the $8$-space $[Z^2 (\mathbb R^{4})]_{U_8}$ consisting of algebraic restrictions of the form $[c_1\theta _1 + \cdots + c_8\theta _8]_{U_8}$ such that $(c_1,c_2,c_3)\ne (0,0,0)$.
\end{theorem}

\begin{theorem}
\label{klasu8} $ \ $

\smallskip

\noindent (i) \ Any algebraic restriction in $[Z ^2(\mathbb R^{2n})]_{U_8}$ can be brought by a symmetry of $U_8$ to one of the normal forms $[U_8]^i$ given in the second column of Table \ref{tabu8}.

\smallskip

\noindent (ii)  \ The codimension in $[Z ^2(\mathbb R^{2n})]_{U_8}$ of the singularity class corresponding to the normal form $[U_8]^i$ is equal to $i$, the symplectic multiplicity and the index of isotropy are given in the fourth and fifth columns of Table  \ref{tabu8}.

\smallskip

\noindent (iii) \ The singularity classes corresponding to the normal forms are disjoint.

\smallskip

\noindent (iv) \ The parameters $c, c_1, c_2$ of the normal forms $[U_8]^i$ are moduli.

\end{theorem}

\renewcommand*{\arraystretch}{1.3}
\begin{center}
\begin{table}[h]
    \begin{small}
    \noindent
    \begin{tabular}{|p{3cm}|p{6.5cm}|c|c|c|}
                      \hline
    symplectic class &   normal forms for algebraic restrictions    & cod & $\mu ^{\rm sym}$ &  ind  \\ \hline
  $(U_8)^0$ \;\;  $(2n\ge 4)$ & $[U_8]^0: [\theta _1 + c_1\theta _2 + c_2\theta _3]_{U_8}$,\;
                          &  $0$ & $2$ & $0$  \\  \hline
  $(U_8)^1$ \;\; $(2n\ge 4)$ & $[U_8]^1: [\pm\theta _2 + c_1\theta _3 + c_2\theta _4]_{U_8}$
                                         &  $1$ & $3$ & $0$ \\ \hline
    $(U_8)^2$ \;\; $(2n\ge 4)$& $[U_8]^2: [\theta _3 + c_1\theta_4+c_2\theta_5]_{U_8}$,
    \; $c_1\ne -\frac{1}{3}, c_1\ne 2$      & $2$ & $4$ & $0$     \\ \hline

    $(U_8)^{3,0}_5$ \;\; $(2n\ge 4)$& $[U_8]^{3,0}_5: [\theta _3 -\frac{1}{3} \theta_4+c_1\theta_5+c_2\theta_6 ]_{U_8}$        & $3$ & $5$ & $0$     \\ \hline

    $(U_8)^{3,0}_\infty$ \;\; $(2n\ge 4)$& $[U_8]^{3,0}_\infty: [\theta _3 + 2\theta_4+c_1\theta_5+c_2\theta_7 ]_{U_8}$
                  & $3$ & $5$ & $0$     \\ \hline
    $(U_8)^{3,1}$ \;\; $(2n\ge 6)$ & $[U_8]^{3,1}: [\theta _4 + c\theta _5]_{U_8}$ &
                $3$ & $4$ & $1$  \\ \hline
    $(U_8)^4$ \;\; $(2n\ge 6)$ & $[U_8]^4: [\pm\theta _5 + c\theta _6]_{U_8}$
                                         &  $4$ & $5$ & $1$    \\ \hline
    $(U_8)^5$ \;\; $(2n\ge 6)$ & $[U_8]^5: [\theta _6 + c\theta _7]_{U_8}$
                                         &  $5$ & $6$ & $2$ \\ \hline
    $(U_8)^6$ \;\; $(2n\ge 6)$ & $[U_8]^6: [\pm\theta _7 + c\theta _8]_{U_8}$
                                         &  $6$ & $7$ & $2$ \\ \hline
    $(U_8)^7$ \;\; $(2n\ge 6)$ & $[U_8]^7: [\theta _8]_{U_8}$
                                         &  $7$ & $7$ & $3$   \\ \hline
    $(U_8)^{8}$ \;\; $(2n\ge 6)$ & $[U_8]^{8}: [0]_{U_8}$ &  $8$ & $8$ & $\infty $ \\ \hline
\end{tabular}

\smallskip

\caption{\small Classification of symplectic $U_8$ singularities.\newline
$cod$ -- codimension of the classes;  $\mu ^{sym}$-- symplectic multiplicity; \newline  $ind$ -- the index of isotropy.}\label{tabu8}
\end{small}
\end{table}
\end{center}


The proof of Theorem \ref{klasu8} is presented in Section \ref{u8-proof}.
\subsubsection{Symplectic normal forms}
\label{u8-normal}

Let us transfer the normal forms $[U_8]^i$  to symplectic normal forms
. Fix a family $\omega ^i$ of symplectic forms on $\mathbb R^{2n}$ realizing the family $[U_8]^i$ of algebraic restrictions.  We can fix, for example,

\smallskip
\begin{small}

\noindent $\omega ^0 = \theta _1 + c_1\theta _2 + c_2\theta _3 +
dx_2\wedge dx_4 + dx_5\wedge dx_6 + \cdots + dx_{2n-1}\wedge
dx_{2n};$

\smallskip

\noindent $\omega ^1 = \pm \theta _2 + c_1\theta _3 + c_2\theta _4 +
dx_1\wedge dx_4 + dx_5\wedge dx_6 + \cdots + dx_{2n-1}\wedge
dx_{2n}; $

\smallskip

\noindent $\omega ^2 = \theta_3+ c_1 \theta _4 + c_2 \theta _5 + dx_3\wedge dx_4 +
dx_5\wedge dx_6 + \cdots + dx_{2n-1}\wedge dx_{2n};$

\smallskip

\noindent $\omega ^{3,0}_5 = \theta_3-\frac{1}{3} \theta _4 + c_1 \theta _5 + c_2 \theta _6+ dx_3\wedge dx_4 + dx_5\wedge dx_6 + \cdots + dx_{2n-1}\wedge dx_{2n};$

\smallskip

\noindent $\omega ^{3,0}_\infty = \theta_3+2 \theta _4 + c_1 \theta _5 + c_2 \theta _7+ dx_3\wedge dx_4 + dx_5\wedge dx_6 + \cdots + dx_{2n-1}\wedge dx_{2n};$

\smallskip

\noindent $\omega ^{3,1} = \theta _4 + c\theta _5 + dx_1\wedge dx_4 +
dx_2\wedge dx_5 + dx_3\wedge dx_6 + dx_7\wedge dx_8 + \cdots +
dx_{2n-1}\wedge dx_{2n};$

\smallskip

\noindent $\omega ^4 = \pm\theta _5 + c\theta _6 + dx_1\wedge dx_4 +dx_2\wedge dx_5 +
dx_3\wedge dx_6 + dx_7\wedge dx_8 + \cdots + dx_{2n-1}\wedge dx_{2n};$

\smallskip

\noindent $\omega ^5 = \theta _6 + c\theta _7+ dx_1\wedge dx_4 + dx_2\wedge dx_5 +
dx_3\wedge dx_6 + dx_7\wedge dx_8 + \cdots + dx_{2n-1}\wedge
dx_{2n};$

\smallskip

\noindent $\omega ^6 = \pm\theta _7 + c\theta _8+ dx_1\wedge dx_4 + dx_2\wedge dx_5 +
dx_3\wedge dx_6 + dx_7\wedge dx_8 + \cdots + dx_{2n-1}\wedge
dx_{2n};$

\smallskip

\noindent $\omega ^7 = \theta _8+ dx_1\wedge dx_4 + dx_2\wedge dx_5 +
dx_3\wedge dx_6 + dx_7\wedge dx_8 + \cdots + dx_{2n-1}\wedge
dx_{2n};$

\smallskip

\noindent $\omega ^8 = dx_1\wedge dx_4 + dx_2\wedge dx_5 +
dx_3\wedge dx_6 + dx_7\wedge dx_8 + \cdots + dx_{2n-1}\wedge
dx_{2n}.$
\end{small}

\smallskip

  Fix, for $i=0,1,\cdots ,8$ a family $\Phi ^i$ of local diffeomorphisms which bring the family of symplectic forms $\omega ^i$ to the symplectic form $\omega_0$: $(\Phi ^i)^*\omega ^i = \omega_0$. Consider the families $U_8^i = (\Phi ^i)^{-1}(U_8)$. Any stratified submanifold of the symplectic space $(\mathbb R^{2n},\omega_0)$ which is diffeomorphic to $U_8$ is symplectically equivalent to one and only one of the normal forms $U_8^i, i = 0,1,\cdots ,8$  presented in Theorem \ref{u8-main}. By Theorem \ref{klasu8} we obtain that  parameters $c,c_1,c_2$ of the normal forms are moduli.


\subsubsection{Proof of Theorem \ref{klasu8}}
\label{u8-proof}

\medskip

\noindent The germ of a vector field tangent to $U_8$ of non trivial action on algebraic restrictions of closed 2-forms to  $U_8$ may be described as a linear combination of germs of vector fields: 
$X_0\!=\!E$,\, $X_1\!=x_3E,\, X_2\!=x_1E,\,  X_3\!=x_3^2E,\, X_4\!=x_2E,\, X_5\!=x_1x_3E,$  $X_6\!=x_3^3E$,  $X_7\!=x_1^2E$, 
$X_8\!=\!x_2x_3E$, where $E$ is the Euler vector field    \begin{equation}
\label{eu8}
E\!=\!3 x_1 \partial /\partial x_1\!+4 x_2 \partial /\partial x_2\!+2 x_3 \partial /\partial x_3.
\end{equation}

\begin{proposition} \label{u8-infinitesimal}

The infinitesimal action of germs of quasi-homogeneous vector
fields tangent to $N\in (U_8)$ on the basis of the vector space of
algebraic restrictions of closed $2$-forms to $N$ is presented in
Table \ref{infini-u8}.

\setlength{\tabcolsep}{1.7mm}
\renewcommand*{\arraystretch}{1.3}
\begin{small}
\begin{table}[h]
\begin{center}
\begin{tabular}{|l|r|r|r|r|r|r|r|r|}

 \hline

  $\mathcal L_{X_i} [\theta_j]$ & $[\theta_1]$   &   $[\theta_2]$ &   $[\theta_3]$ & $[\theta_4]$   & $[\theta_5]$   & $[\theta_6]$   & $[\theta_7]$ & $[\theta_8]$ \\ \hline

  $X_0=E$ & $5 [\theta_1]$ & $6 [\theta_2]$ & $7 [\theta_3]$ & $7 [\theta_4]$ & $8 [\theta_5]$ & $9[\theta_6]$ & $10 [\theta_7]$ & $11 [\theta_8]$ \\ \hline

  $X_1=x_3E$ & $7[\theta_4]$ & $-16 [\theta_5]$ & $ 3[\theta_6]$ & $9 [\theta_6]$ & $10[\theta_7]$  & $11[\theta_8]$ & $[0]$ & $[0]$ \\  \hline

  $X_2=x_1E$ & $8[\theta_5]$  & $ -6[\theta_6]$ & $ -20[\theta_7]$ & $10[\theta_7]$ & $\frac{11}{3}[\theta_8]$ & $[0]$ & $[0]$ & $[0]$ \\  \hline

  $X_3=x_3^2E$ & $9[\theta_6]$ & $-20[\theta_7]$  & $\frac{11}{3}[\theta_8]$ & $11[\theta_8]$ & $[0]$ & $[0]$ & $[0]$ & $[0]$ \\ \hline

  $X_4=x_2 E$ & $ -3[\theta_6]$  & $-40[\theta_7]$ & $-\frac{55}{3}[\theta_8]$  & $-\frac{11}{3}[\theta_8]$ & $[0]$ & $[0]$  & $[0]$ & $[0]$ \\  \hline

  $X_5=x_1x_3 E$ & $10[\theta_7]$  & $-\frac{22}{3}[\theta_8]$ & $[0]$  & $[0]$ & $[0]$   & $[0]$ & $[0]$ & $[0]$ \\  \hline

  $X_6=x_3^3 E$ & $11[\theta_8]$  & $[0]$ & $[0]$  & $[0]$ & $[0]$   & $[0]$ & $[0]$ & $[0]$ \\  \hline

  $X_7=x_1^2 E$ & $\frac{11}{3}[\theta_8]$  & $[0]$ & $[0]$  & $[0]$ & $[0]$   & $[0]$ & $[0]$ & $[0]$ \\  \hline

  $X_8=x_2x_3 E$ & $-\frac{11}{3}[\theta_8]$  & $[0]$ & $[0]$  & $[0]$ & $[0]$   & $[0]$ & $[0]$ & $[0]$ \\  \hline

\end{tabular}
\end{center}

\caption{\small Infinitesimal actions on algebraic restrictions of closed \newline
2-forms to   $U_8$. ($E$ is defined as in (\ref{eu8}).)}\label{infini-u8}
\end{table}
\end{small}
\end{proposition}


\medskip

Let $\mathcal{A}=[c_1 \theta_1+c_2 \theta_2+c_3 \theta_3+c_4 \theta_4+c_5 \theta_5+c_6 \theta_6 +c_7 \theta_7+c_8 \theta_8]_{U_8}$
be the algebraic restriction of a symplectic form $\omega$.

\medskip

The first statement of Theorem \ref{klasu8} follows from the following lemmas.

\begin{lemma}
\label{u8lem0} If \;$c_1\ne 0$\; then the algebraic restriction
$\mathcal{A}=[\sum_{k=1}^8 c_k \theta_k]_{U_8}$
can be reduced by a symmetry of $U_8$ to an algebraic restriction $[\theta_1+\widetilde{c}_2 \theta_2+\widetilde{c}_3 \theta_3]_{U_8}$.
\end{lemma}




\begin{lemma}
\label{u8lem1} If \;$c_1\!=\!0$\;and $c_2\ne\!0$ then the algebraic restriction  $\mathcal{A}$ 
can be reduced by a symmetry of $U_8$ to an algebraic restriction $[\pm\theta_2+\widetilde{c}_3 \theta_3+\widetilde{c}_4 \theta_4]_{U_8}$.
\end{lemma}

\begin{lemma}
\label{w8lem2} If $c_1=c_2=0$ and $c_3\ne0$, $c_4\ne 2c_3$, $c_4\ne -\frac{1}{3}c_3$ then the algebraic restriction $\mathcal{A}$
can be reduced by a symmetry of \, $U_8$ to an algebraic restriction $[\theta_3+\widetilde{c}_4 \theta_4+\widetilde{c}_5 \theta_5]_{U_8}$.
\end{lemma}

\begin{lemma}
\label{w8lem3a} If $c_1=c_2=0$ and $c_3\ne0$, $c_4= -\frac{1}{3}c_3$ then the algebraic restriction $\mathcal{A}$
can be reduced by a symmetry of \, $U_8$ to an algebraic restriction $[\theta_3-\frac{1}{3} \theta_4+\widetilde{c}_5 \theta_5+\widetilde{c}_6 \theta_6]_{U_8}$.
\end{lemma}

\begin{lemma}
\label{w8lem3b} If $c_1=c_2=0$ and $c_3\ne0$, $c_4= 2c_3$ then the algebraic restriction $\mathcal{A}$
can be reduced by a symmetry of \, $U_8$ to an algebraic restriction $[\theta_3+2\theta_4+\widetilde{c}_5 \theta_5+\widetilde{c}_7 \theta_7]_{U_8}$.
\end{lemma}

\begin{lemma}
\label{u8lem31} If \;$c_1=c_2=c_3=0$ and $c_4\ne 0$\; then the algebraic restriction $\mathcal{A}$
can be reduced by a symmetry of \, $U_8$ to an algebraic restriction   $[\theta_4+\widetilde{c}_5 \theta_5]_{U_8}$.
\end{lemma}

\begin{lemma}
\label{u8lem4} If \;$c_1=0,\ldots,c_4=0$ and $c_5\ne 0$,\; then the algebraic restriction $\mathcal{A}$
can be reduced by a symmetry of \, $U_8$ to an algebraic restriction $[\pm\theta_5+\widetilde{c}_6 \theta_6]_{U_8}$.
\end{lemma}

\begin{lemma}
\label{u8lem5} If \;$c_1=0,\ldots,c_5=0$ and $c_6\ne 0$\; then the algebraic restriction $\mathcal{A}$
can be reduced by a symmetry of \, $U_8$ to an algebraic restriction $[\theta_6+\widetilde{c}_7 \theta_7]_{U_8}$.
\end{lemma}

\begin{lemma}
\label{u8lem6} If \;$c_1=0,\ldots,c_6=0$ and $c_7\ne 0$\; then the algebraic restriction $\mathcal{A}$
can be reduced by a symmetry of \, $U_8$ to an algebraic restriction $[\pm\theta_7+\widetilde{c}_8 \theta_8]_{U_8}$.
\end{lemma}

\begin{lemma}
\label{u8lem7} If \;$c_1=0,\ldots,c_7=0$ and $c_8\ne 0$\; then the algebraic restriction $\mathcal{A}$
can be reduced by a symmetry of \, $U_8$ to an algebraic restriction $[\theta_8]_{U_8}$.
\end{lemma}

The proofs of Lemmas \ref{u8lem0} -- \ref{u8lem7} are similar and based on Table \ref{infini-u7}, Proposition \ref{elimin1} or the homotopy method.

\medskip

To prove statement $(iii)$ of Theorem \ref{klasu8} we have to show that singularity classes corresponding to normal forms are disjoint. The singularity classes that can be distinguished by geometric conditions obviously are disjoint. From Theorem \ref{geom-cond-u8} we see that only classes $(U_8)^2$ and $(U_8)^{3,0}_5$ can not be distinguished by the geometric conditions but their symplectic multiplicities are distinct, hence the classes are disjoint.

The proofs of statements $(ii)$ and $(iv)$ of Theorem \ref{klasu8} are similar to analogous proofs for Theorem \ref{klasu7}. 



\pagebreak

\subsection{Proofs for  $U_9$ singularity}
\subsubsection{Algebraic restrictions to $U_9$ and their classification}\label{u9-class} $ \ $

One has the following relations for $(U_9)$-singularities
\begin{equation}
[x_1^2+x_2x_3]_{U_9}=0.
\label{u901}
\end{equation}
\begin{equation}
[x_1x_2+x_3^4]_{U_9}=0,
\label{u902}
\end{equation}
\begin{equation}
[d(x_1^2+x_2x_3)]_{U_9}=[2x_1dx_1+ x_2dx_3+x_3 dx_2]_{U_9}=0
\label{u91}
\end{equation}
\begin{equation}
[d(x_1x_2+x_3^4)]_{U_9}=[x_1dx_2+x_2dx_1+4x_3^3dx_3]_{U_9}=0
\label{u92}
\end{equation}
Multiplying these relations by suitable $1$-forms and $2$-forms we obtain the relations towards calculating $[\Lambda^2(\mathbb R^{2n})]_N$ for $N=U_9$.

\begin{proposition}
\label{u9-all}
The space $[\Lambda ^{2}(\mathbb R^{2n})]_{U_9}$ is a $10$-dimensional vector space spanned by the algebraic restrictions to $U_9$ of the $2$-forms

$\theta _1= dx_1\wedge dx_3, \;\; \theta _2=dx_2\wedge dx_3,\;\; \theta_3 = dx_1\wedge dx_2,$


  $\theta _4 = x_3dx_1\wedge dx_3,\;\; \theta _5 = x_1dx_1\wedge dx_3,$\;\; $\theta _6= x_3^2 dx_1\wedge dx_3$, \; \; $\sigma = x_3 dx_1\wedge dx_2,$


 $\theta _7= x_1x_3 dx_1\wedge dx_3$,\;\; $\theta _8= x_3^3 dx_1\wedge dx_3$,\;\; $\theta _9= x_1x_3^2 dx_1\wedge dx_3$.
\end{proposition}

Proposition \ref{u9-all} and results of Section \ref{method}  imply the following description of the space $[Z ^2(\mathbb R^{2n})]_{U_9}$ and the manifold $[{\rm Symp} (\mathbb R^{2n})]_{U_9}$.

\begin{theorem} \label{u9-baza}
 The space $[Z^2(\mathbb R^{2n})]_{U_9}$ is a $9$-dimensional vector space
 spanned by the algebraic restrictions to $U_9$  of the quasi-homogeneous $2$-forms $\theta_i$  of degree $\delta_i$

\begin{small}

$\theta _1= dx_1\wedge dx_3,\;\;\;\delta_1=8,$

\smallskip

$\theta _2=dx_2\wedge dx_3,\;\;\;\delta_2=10,$

\smallskip

$\theta_3 = dx_1\wedge dx_2,\;\;\;\delta_3=12,$

\smallskip

  $\theta _4 = x_3dx_1\wedge dx_3,\;\;\;\delta_4=11,$

  \smallskip

  $\theta _5 = x_1dx_1\wedge dx_3,\;\;\;\delta_5=13,$

\smallskip

$\theta _6=x_3^2 dx_1\wedge dx_3,\;\;\;\delta_6=14, $

\smallskip

$\theta _7= x_1x_3 dx_1\wedge dx_3,\;\;\;\delta_7=16$,

\smallskip

$\theta _8= x_3^3 dx_1\wedge dx_3,\;\;\;\delta_8=17$,

\smallskip

$\theta _9= x_1x_3^2 dx_1\wedge dx_3,\;\;\;\delta_9=19$,
\end{small}

\smallskip

If $n\ge 3$ then $[{\rm Symp} (\mathbb R^{2n})]_{U_9} = [Z^2(\mathbb R^{2n})]_{U_9}$. The manifold $[{\rm Symp} (\mathbb R^{4})]_{U_9}$ is an open part of the $9$-space $[Z^2 (\mathbb R^{4})]_{U_9}$ consisting of algebraic restrictions of the form $[c_1\theta _1 + \cdots + c_9\theta _9]_{U_9}$ such that $(c_1,c_2,c_3)\ne (0,0,0)$.
\end{theorem}

\begin{theorem}
\label{klasu9} $ \ $

\smallskip

\noindent (i) \ Any algebraic restriction in $[Z ^2(\mathbb R^{2n})]_{U_9}$ can be brought by a symmetry of $U_9$ to one of the normal forms $[U_9]^i$ given in the second column of Table \ref{tabu9}.

\smallskip

\noindent (ii)  \ The codimension in $[Z ^2(\mathbb R^{2n})]_{U_9}$ of the singularity class corresponding to the normal form $[U_9]^i$ is equal to $i$, the symplectic multiplicity and the index of isotropy are given in the fourth and fifth columns of Table  \ref{tabu9}.

\smallskip

\noindent (iii) \ The singularity classes corresponding to the normal forms are disjoint.

\smallskip

\noindent (iv) \ The parameters $c, c_1, c_2, c_3$ of the normal forms $[U_9]^i$ are moduli.

\end{theorem}

\renewcommand*{\arraystretch}{1.3}
\begin{center}
\begin{table}[h]
    \begin{small}
    \noindent
    \begin{tabular}{|p{3cm}|p{6.5cm}|c|c|c|}
                      \hline
    symplectic class &   normal forms for algebraic restrictions    & cod & $\mu ^{\rm sym}$ &  ind  \\ \hline
  $(U_9)^0$ \;\;  $(2n\ge 4)$ & $[U_9]^0: [\pm\theta _1 + c_1\theta _2 + c_2\theta _3]_{U_9}$,\;
                          &  $0$ & $2$ & $0$  \\  \hline
  $(U_9)^1$ \;\; $(2n\ge 4)$ & $[U_9]^1: [\pm\theta _2 + c_1\theta _3 + c_2\theta _4+c_3\theta_6]_{U_9}$
                                         &  $1$ & $4$ & $0$ \\ \hline
    $(U_9)^2$ \;\; $(2n\ge 4)$& $[U_9]^2: [\pm\theta _3 + c_1\theta_4+c_2\theta_5]_{U_9}$,
    \; $c_1\ne0$      & $2$ & $4$ & $0$     \\ \hline

    $(U_9)^{3,0}$ \;\; $(2n\ge 4)$& $[U_9]^{3,0}: [\pm\theta _3 +c_1\theta_5+c_2\theta_6 ]_{U_9}$,
    \; $c_1\ne0$           & $3$ & $5$ & $0$     \\ \hline

    $(U_9)^{4,0}$ \;\; $(2n\ge 4)$& $[U_9]^{4,0}: [\pm\theta _3 + c_1\theta_6+c_2\theta_7 ]_{U_9}$
                  & $4$ & $6$ & $0$     \\ \hline
    $(U_9)^{3,1}$ \;\; $(2n\ge 6)$ & $[U_9]^{3,1}: [\theta _4 + c\theta _5]_{U_9}$ &
                $3$ & $4$ & $1$  \\ \hline
    $(U_9)^{4,1}$ \;\; $(2n\ge 6)$ & $[U_9]^{4,1}: [\theta _5 + c_1\theta _6+c_2\theta _8]_{U_9}$
                                         &  $4$ & $6$ & $1$    \\ \hline
    $(U_9)^5$ \;\; $(2n\ge 6)$ & $[U_9]^5: [\pm\theta _6 + c\theta _7]_{U_9}$
                                         &  $5$ & $6$ & $2$ \\ \hline
    $(U_9)^6$ \;\; $(2n\ge 6)$ & $[U_9]^6: [\pm\theta _7 + c\theta _8]_{U_9}$
                                         &  $6$ & $7$ & $2$ \\ \hline
    $(U_9)^7$ \;\; $(2n\ge 6)$ & $[U_9]^7: [\theta _8+c\theta_9]_{U_9}$
                                         &  $7$ & $8$ & $3$   \\ \hline
    $(U_9)^{8}$ \;\; $(2n\ge 6)$ & $[U_9]^{8}: [\theta_9]_{U_9}$
                                         &  $8$ & $8$ & $3$   \\ \hline
    $(U_9)^{9}$ \;\; $(2n\ge 6)$ & $[U_9]^{9}: [0]_{U_9}$ &  $9$ & $9$ & $\infty $ \\ \hline
\end{tabular}
\smallskip

\caption{\small Classification of symplectic $U_9$ singularities.\newline
$cod$ -- codimension of the classes;  $\mu ^{sym}$-- symplectic multiplicity;  \newline $ind$ -- the index of isotropy.}\label{tabu9}
\end{small}
\end{table}
\end{center}


The proof of Theorem \ref{klasu9} is presented in Section \ref{u9-proof}.

\subsubsection{Symplectic normal forms}
\label{u9-normal} $ \ $

\medskip
\par
Let us transfer the normal forms $[U_9]^i$  to symplectic normal forms.  We fix a family $\omega ^i$ of symplectic forms on $\mathbb R^{2n}$ realizing the family $[U_9]^i$ of algebraic restrictions.

\smallskip
\begin{small}

\noindent $\omega ^0 = \pm\theta _1 + c_1\theta _2 + c_2\theta _3 +
dx_2\wedge dx_4 + dx_5\wedge dx_6 + \cdots + dx_{2n-1}\wedge
dx_{2n};$

\smallskip

\noindent $\omega ^1 = \pm \theta _2 + c_1\theta _3 + c_2\theta _4 +c_3\theta _6 +
dx_1\wedge dx_4 + dx_5\wedge dx_6 + \cdots + dx_{2n-1}\wedge
dx_{2n}; $

\smallskip

\noindent $\omega ^2 = \pm \theta_3+ c_1 \theta _4 + c_2 \theta _5 + dx_3\wedge dx_4 +
dx_5\wedge dx_6 + \cdots + dx_{2n-1}\wedge dx_{2n}, \; c_1\ne 0;$

\smallskip

\noindent $\omega ^{3,0} = \pm\theta_3 + c_1 \theta _5 + c_2 \theta _6+ dx_3\wedge dx_4 + dx_5\wedge dx_6 + \cdots + dx_{2n-1}\wedge dx_{2n}, \; c_1\ne 0;$

\smallskip

\noindent $\omega ^{4,0} = \pm\theta_3+ c_1 \theta _6 + c_2 \theta _7+ dx_3\wedge dx_4 + dx_5\wedge dx_6 + \cdots + dx_{2n-1}\wedge dx_{2n};$

\smallskip

\noindent $\omega ^{3,1} = \theta _4 + c\theta _5 + dx_1\wedge dx_4 +
dx_2\wedge dx_5 + dx_3\wedge dx_6 + dx_7\wedge dx_8 + \cdots +
dx_{2n-1}\wedge dx_{2n};$

\smallskip

\noindent $\omega ^{4,1} = \theta _5 + c_1\theta _6 + c_2\theta _8+ dx_1\wedge dx_4 +dx_2\wedge dx_5 +
dx_3\wedge dx_6 + dx_7\wedge dx_8 + \cdots + dx_{2n-1}\wedge dx_{2n};$

\smallskip

\noindent $\omega ^5 = \pm\theta _6 + c\theta _7+ dx_1\wedge dx_4 + dx_2\wedge dx_5 +
dx_3\wedge dx_6 + dx_7\wedge dx_8 + \cdots + dx_{2n-1}\wedge
dx_{2n};$

\smallskip

\noindent $\omega ^6 = \pm\theta _7 + c\theta _8+ dx_1\wedge dx_4 + dx_2\wedge dx_5 +
dx_3\wedge dx_6 + dx_7\wedge dx_8 + \cdots + dx_{2n-1}\wedge
dx_{2n};$

\smallskip

\noindent $\omega ^7 = \theta _8+ c\theta _9+ dx_1\wedge dx_4 + dx_2\wedge dx_5 +
dx_3\wedge dx_6 + dx_7\wedge dx_8 + \cdots + dx_{2n-1}\wedge
dx_{2n};$

\smallskip

\noindent $\omega ^{8} = \theta _9+ dx_1\wedge dx_4 + dx_2\wedge dx_5 +
dx_3\wedge dx_6 + dx_7\wedge dx_8 + \cdots + dx_{2n-1}\wedge
dx_{2n};$


\noindent $\omega ^{9} = dx_1\wedge dx_4 + dx_2\wedge dx_5 +
dx_3\wedge dx_6 + dx_7\wedge dx_8 + \cdots + dx_{2n-1}\wedge
dx_{2n}.$
\end{small}

\medskip

\subsubsection{Proof of Theorem \ref{klasu9}}
\label{u9-proof}

\medskip

\noindent The germ of a vector field tangent to $U_9$ of non trivial action on algebraic restrictions of closed 2-forms to  $U_9$ may be described as a linear combination of germs of vector fields: 
$X_0\!=\!E, X_1\!=x_3E, X_2\!=x_1E, X_3\!=x_3^2E, X_4\!=x_2E, X_5\!=x_1x_3E,$ $X_6\!=x_3^3E,  X_7\!=x_1x_3^2E$, where $E$ is the Euler vector field    \begin{equation}
\label{eu9}
E=5 x_1 \partial /\partial x_1+7 x_2 \partial /\partial x_2+3 x_3 \partial /\partial x_3.
\end{equation}

\begin{proposition} \label{u9-infinitesimal}

The infinitesimal action of germs of quasi-homogeneous vector
fields tangent to $N\in (U_9)$ on the basis of the vector space of
algebraic restrictions of closed $2$-forms to $N$ is presented in
Table \ref{infini-u9}.
\setlength{\tabcolsep}{1.4mm}
\renewcommand*{\arraystretch}{1.3}
\begin{small}
\begin{table}[h]
\begin{center}
\begin{tabular}{|l|r|r|r|r|r|r|r|r|r|}

 \hline
  $\mathcal L_{X_i} [\theta_j]$ & $[\theta_1]$   &   $[\theta_2]$ &   $[\theta_3]$ & $[\theta_4]$   & $[\theta_5]$   & $[\theta_6]$   & $[\theta_7]$ & $[\theta_8]$ & $[\theta_9]$ \\ \hline

  $X_0\!=\!E$ & $8 [\theta_1]$ & $10 [\theta_2]$ & $12 [\theta_3]$ & $11 [\theta_4]$ & $13 [\theta_5]$ & $14 [\theta_6]$ & $16 [\theta_7]$  & $17 [\theta_8]$ & $19[\theta_9]$ \\ \hline

  $X_1\!=\!x_3E$  & $11[\theta_4]$ & $-26 [\theta_5]$ & $[0]$  & $14[\theta_6]$  & $16[\theta_7]$ & $17[\theta_8]$  & $19[\theta_9]$ & $[0]$ & $[0]$ \\  \hline

  $X_2\!=\!x_1E$ & $13 [\theta_5]$ & $ [0]$ & $-68 [\theta_8]$ & $16 [\theta_7]$ & $ [0]$ & $19[\theta_9]$ & $0[0]$ & $[0]$ & $[0]$ \\  \hline

  $X_3\!=\!x_3^2 E$ & $14 [\theta_6]$  & $-32[\theta_7]$ & $[0]$  & $17[\theta_8]$ & $19[\theta_9]$ & $[0]$  & $[0]$  & $[0]$ & $[0]$ \\  \hline

  $X_4\!=\!x_2E$ & $ [0]$ & $-136[\theta_8]$  & $-38[\theta_9]$ & $[0]$ & $[0]$ & $[0]$ & $[0]$ & $[0]$ & $[0]$ \\ \hline

  $X_5\!=\!x_1x_3 E$ & $16[\theta_7]$  & $[0]$ & $[0]$  & $19[\theta_9]$ & $[0]$ & $[0]$  & $[0]$  & $[0]$ & $[0]$ \\  \hline

  $X_6\!=\!x_3^3 E$ & $17[\theta_8]$  & $-38[\theta_9]$  & $[0]$ & $[0]$ & $[0]$ & $[0]$ & $[0]$ & $[0]$ & $[0]$ \\ \hline

  $X_7\!=\!x_1x_3^2 E$ & $19[\theta_9]$  & $[0]$ & $[0]$  & $[0]$ & $[0]$ & $[0]$  & $[0]$  & $[0]$ & $[0]$ \\  \hline

\end{tabular}
\end{center}

\smallskip

\caption{\small Infinitesimal actions on algebraic restrictions of closed\newline
2-forms to   $U_9$. ($E$ is defined as in (\ref{eu9}).)}\label{infini-u9}
\end{table}
\end{small}
\end{proposition}

\medskip

Let $\mathcal{A}=[c_1 \theta_1+c_2 \theta_2+c_3 \theta_3+c_4 \theta_4+c_5 \theta_5+c_6 \theta_6 +c_7 \theta_7+c_8 \theta_8+c_9 \theta_9]_{U_9}$
be the algebraic restriction of a symplectic form $\omega$.

The first statement of Theorem \ref{klasu9} follows from the following lemmas.

\begin{lemma}
\label{u9lem0} If \;$c_1\ne 0$\; then the algebraic restriction
$\mathcal{A}=[\sum_{k=1}^9 c_k \theta_k]_{U_9}$
can be reduced by a symmetry of $U_9$ to an algebraic restriction $[\pm\theta_1+\widetilde{c}_2 \theta_2+\widetilde{c}_3 \theta_3]_{U_9}$.
\end{lemma}

\begin{lemma}
\label{u9lem1} If \;$c_1\!=\!0$\;and $c_2\ne\!0$ then the algebraic restriction  $\mathcal{A}$ 
can be reduced by a symmetry of $U_9$ to an algebraic restriction $[\pm\theta_2+\widetilde{c}_3 \theta_3+\widetilde{c}_4 \theta_4+\widetilde{c}_6\theta_6]_{U_9}$.
\end{lemma}

\begin{lemma}
\label{w9lem2} If $c_1=c_2=0$ and $c_3\cdot c_4\ne 0$ then the algebraic restriction $\mathcal{A}$
can be reduced by a symmetry of \, $U_9$ to an algebraic restriction $[\pm\theta_3+\widetilde{c}_4 \theta_4+ \widetilde{c}_5\theta_5]_{U_9}$.
\end{lemma}

\begin{lemma}
\label{w9lem30} If $c_1=c_2=c_4=0$ and $c_3\cdot c_5\ne 0$ then the algebraic restriction $\mathcal{A}$
can be reduced by a symmetry of \, $U_9$ to an algebraic restriction $[\pm\theta_3+\widetilde{c}_5 \theta_5+ \widetilde{c}_6\theta_6]_{U_9}$.
\end{lemma}

\begin{lemma}
\label{w9lem40} If $c_1=c_2=c_4=c_5=0$ and $c_3\ne 0$ then the algebraic restriction $\mathcal{A}$
can be reduced by a symmetry of \, $U_9$ to an algebraic restriction $[\pm\theta_3+\widetilde{c}_6 \theta_6+ \widetilde{c}_7\theta_7]_{U_9}$.
\end{lemma}

\begin{lemma}
\label{u9lem31} If \;$c_1=c_2=c_3=0$ and $c_4\ne 0$\; then the algebraic restriction $\mathcal{A}$
can be reduced by a symmetry of \, $U_9$ to an algebraic restriction   $[\theta_4+\widetilde{c}_5 \theta_5]_{U_9}$.
\end{lemma}

\begin{lemma}
\label{u8lem41} If \;$c_1=0,\ldots,c_4=0$ and $c_5\ne 0$,\; then the algebraic restriction $\mathcal{A}$
can be reduced by a symmetry of \, $U_9$ to an algebraic restriction $[\theta_5+\widetilde{c}_6 \theta_6+\widetilde{c}_8 \theta_8]_{U_9}$.
\end{lemma}

\begin{lemma}
\label{u9lem5} If \;$c_1=0,\ldots,c_5=0$ and $c_6\ne 0$\; then the algebraic restriction $\mathcal{A}$
can be reduced by a symmetry of \, $U_9$ to an algebraic restriction $[\pm\theta_6+\widetilde{c}_7 \theta_7]_{U_9}$.
\end{lemma}

\begin{lemma}
\label{u9lem6} If \;$c_1=0,\ldots,c_6=0$ and $c_7\ne 0$\; then the algebraic restriction $\mathcal{A}$
can be reduced by a symmetry of \, $U_9$ to an algebraic restriction $[\pm\theta_7+\widetilde{c}_8 \theta_8]_{U_9}$.
\end{lemma}

\begin{lemma}
\label{u9lem7} If \;$c_1=0,\ldots,c_7=0$ and $c_8\ne 0$\; then the algebraic restriction $\mathcal{A}$
can be reduced by a symmetry of \, $U_9$ to an algebraic restriction $[\theta_8+\widetilde{c}_9 \theta_9]_{U_9}$.
\end{lemma}

\begin{lemma}
\label{u9lem8} If \;$c_1=0,\ldots,c_8=0$ and $c_9\ne 0$\; then the algebraic restriction $\mathcal{A}$
can be reduced by a symmetry of \, $U_9$ to an algebraic restriction $[\theta_9]_{U_9}$.
\end{lemma}

The proofs of Lemmas \ref{u9lem0} -- \ref{u9lem8} are similar and based on Table \ref{infini-u9}, Proposition \ref{elimin1} or the homotopy method.

\medskip

The proofs of statements $(ii)$ -- $(iv)$ of Theorem \ref{klasu9} are similar to analogous proofs for Theorem \ref{klasu7}.

\bigskip



\end{document}